\documentclass[11pt,a4paper]{article}
\usepackage{fontenc,amsmath,mathrsfs,amsfonts,amsthm, verbatim,geroENG,amssymb, MnSymbol, enumerate, multicol, color, hyperref, graphicx}
\usepackage{pgf,tikz}
\usepackage{mathrsfs}
\usepackage[english]{babel}
\usepackage[top=1.25in, bottom=1.25in, left=1.0in, right=1.0in]{geometry}
\makeatletter
\@namedef{subjclassname@2010}{%
  \textup{2010} Mathematics Subject Classification}
\makeatother

\begin{document}

\title{Purely Periodic and Transcendental Complex Continued Fractions}

\author{Gerardo Gonz\'alez Robert}

\date{}

\maketitle
\begin{abstract}
Adolf Hurwitz proposed in 1887 a continued fraction algorithm for complex numbers: Hurwitz continued fractions (HCF). Among other similarities between HCF and regular continued fractions, quadratic irrational numbers over $\QU(i)$ are precisely those with periodic HCF expansions (\cite{hur87}, p.196). In this paper, we give some necessary as well as some sufficient conditions for pure periodicity of HCF. Then, we characterize badly approximable complex numbers in terms of HCF. Finally, we prove a slightly weaker complex analogue of a theorem by Y. Bugeaud (\cite{bug13autom}, Theorem 3.1.) on the transcendence of certain continued fractions.
\end{abstract}



\section{Introduction}

Regular continued fractions are a remarkably useful tool in number theory. The structure of a regular continued fraction sometimes helps to determine algebraic or analytic properties of the number it represents. A famous result in this direction is the Euler--Lagrange Theorem, which states that an irrational number has a periodic continued fraction if and only if it is a quadratic surd. Another important theorem, due to Liouville, allows us to construct transcendental numbers by simply taking continued fractions whose terms grow fast enough (\cite{khin}, Theorem 27). A well known conjecture relating the boundedness of the partial quotients with the transcendence of the limit is the following.
\begin{conj01}[Folklore Conjecture]
If the regular continued fraction of an algebraic irrational number $x$ is bounded, then $x$ is a quadratic irrational.
\end{conj01}
Although the conjecture remains widely open, there are important partial results. On the basis of Roth's Theorem, Alan Baker showed in \cite{baker} that whenever the continued fraction of a number $x\in\RE\setminus\QU$ satisfies certain combinatorial condition, $x$ is transcendental. In \cite{quef98}, Martine Qu\'effelec showed the transcendence of the numbers whose regular continued fraction is the Thue-Morse sequence over any alphabet $\{a,b\}\subseteq\Na$. Afterwards, she showed in \cite{quef00} the transcendence of a larger class of automatic continued fractions. Based on their joint work with Florian Luca on $b$-ary expansions in \cite{adbulu04}, Boris Adamczewski and Yann Bugeaud generalized in \cite{adbug05} Qu\'effelec's work by providing weaker sufficient conditions for the transcendence of continued fractions. These results were crowned by Y. Bugeaud in \cite{bug13autom}. The main result of the later paper implies that any real number with an automatic continued fraction is either a  quadratic irrational or transcendental.

Recently, Y. Bugeaud and Dong Han Kim defined in \cite{bugkim} a function giving another notion of complexity of an infinite word. As they noted, their function makes it possible to state the main theorem of \cite{bug13autom} (Theorem \ref{BugeaudRep} below) in an extremely neat fashion. 

For $j,k\in\Za$ satsifying $j\leq k$, we write $[j..k]=[j,k]\cap\Za$.

\begin{def01}\label{DefRepExp}
Let $\clA\neq \vac$ be a finite set and $\bfa=a_1a_2a_3\ldots$ an infinite word on $\clA$. The \textbf{repetition function}, $r(\cdot,\bfa):\Na\to\Na$, is given by
\[
r(n,\bfa) = \min\left\{ m\in\Na_{\geq n+1}: \exists i\in[1..m-n] \; a_{i}\cdots a_{i+n-1}=a_{m-n+1}\cdots a_m\right\}
\]
for all $n\in\Na$. The \textbf{repetition exponent} of $\bfa$, $\rep(\bfa)$, is
\[
\rep(\bfa):= \liminf_{n\to\infty} \frac{r(n,\bfa)}{n}.
\]
\end{def01}

\begin{teo01}[Y. Bugeaud, \cite{bug13autom}]\label{BugeaudRep}
Let $\bfa=a_1a_2\ldots$ be an non-periodic infinite word on a finite subset of $\Na$. If 
\[
\rep( \bfa)< + \infty,
\]
then $\alpha=[0;a_1,a_2,a_3,\ldots]$ is transcendental.
\end{teo01}
\begin{rema1000}
The original result, Theorem 3.1. of \cite{bug13autom}, allows $\bfa$ to take infinitely many values. However, an additional restriction on the corresponding sequence of continuants $\seqn$ is needed, namely $\sup_{n\in\Na} q_{n}^{1/n}<+\infty$. This condition holds when $\bfa$ is bounded. 
\end{rema1000}

An automatic sequence (or $k$-automatic sequence) is a sequence generated by a finite automaton (for a precise definition see \cite{alshal}, Definition 5.1.1.). By Cobham's Theorem on the complexity of automatic sequences (\cite{alshal}, Corollary 10.3.2) and Lemma 2.2. in \cite{bugkim}, every automatic sequence has a finite repetition exponent. Hence, Theorem \ref{BugeaudRep} gives the next result.

\begin{teo01}[Y. Bugeaud, \cite{bug13autom}]\label{BugeaudAutom}
Let $(a_n)_{n=1}^{\infty}$ be an automatic sequence in a finite subset of $\Na$. If $(a_n)_{n=1}^{\infty}$ is not periodic, then $\alpha=[0;a_1,a_2,\ldots]$ is transcendental.
\end{teo01}

It is natural to look for analogues of regular continued fractions in other contexts. An outstanding example for the complex plane $\Cx$ was given by Asmus Schmidt in \cite{aschm01}. His sophisticated construction focuses on the quality of approximation. A much simpler algorithm was proposed by Adolf Hurwitz in \cite{hur87}, it is just the straightforward generalization of a nearest integer continued fraction algorithm (see Section \ref{SecHCF} for details). Some other expansions were studied by Julius Hurwitz \cite{juhu}, William LeVeque \cite{leveque}, Shigeru Tanaka \cite{tanaka}, Georges Poitou \cite{poitou}, among others. Lately, some families of complex continued fractions have been studied by Shrikrishna Gopalrao Dani and Arnaldo Nogueira in \cite{daninog}, \cite{dani}. In this paper, we restrict ourselves to Hurwitz continued fractions. 

Although similarities between Hurwitz and regular continued fractions abound, there are important differences. For instance, Serret's theorem states that two real numbers belong to the same orbit of $\PSL(2,\Za)/\RE$ (acting via M\"obius transformations) if and only if they are both rational or if they are both irrational and their continued fraction expansion eventually coincide. The analogue fails in $\Cx$ with Hurwitz continued fractions. In \cite{lakein03}, Richard Lakein gave a pair of $\PGL(2,\Za[i])$-equivalent complex numbers such that the tails of their Hurwitz continued fraction never coincide.

The most striking result was obtained by Doug Hensley in \cite{hensley} and extended by Wieb Bosma and David Gruenewald in \cite{bosgruen}. It is a negative answer to the Folklore Conjecture for complex numbers and Hurwitz continued fractions.

\begin{teo01}[W. Bosma, D. Gruenewald, \cite{bosgruen}]\label{TeoBoGru}
Let $n$ be a natural number. There exists an algebraic number $\alpha\in\Cx$ whose Hurwitz continued fraction has bounded partial quotients and such that $\alpha$ is of degree $2n$ over $\QU(i)$; that is
\[
[\QU(i,\alpha):\QU(i)]=2n.
\]
\end{teo01}

In spite of the Theorem \ref{TeoBoGru}, we can conclude certain properties about the repetition exponent of the continued fractions of algebraic numbers. A trivial consequence of our main result, Theorem \ref{TeoGero01}, is a weaker complex version of Theorem \ref{BugeaudRep} (see below for notation).

\begin{teo01}\label{TeoGero}
Let $\bfa=a_1a_2a_3\ldots$ be a non-periodic infinite word on $\Za[i]$. If 
\[
\sqrt{8}\leq \liminf_{n\to\infty} |a_n|, \quad
\rep( \bfa)< + \infty,
\]
then $\zeta=[0;a_1,a_2,\ldots]$ is transcendental.
\end{teo01}

\begin{coro01}\label{CoroGero}
Let $\bfa=a_1a_2a_3\ldots$ be an infinite word on $\Za[i]$. If $|a_n|\geq \sqrt{8}$ for all $n\in\Na$ and $\bfa$ is automatic, then $\zeta = [0;a_1,a_2,\ldots]$ is quadratic over $\QU(i)$ or transcendental.
\end{coro01}

The paper is organized as follows. Section 2 defines Hurwitz continued fractions and discusses the associated shift space. Section 3 gives some necessary as well as some sufficient conditions for pure periodicity of Hurwitz continued fractions (Theorem \ref{PurelyPeriodicHCF}). Section 4 characterizes badly approximable complex numbers in terms of Hurwitz continued fractions (Theorem \ref{BadCCharac}). This result was recently shown by Robert Hines in \cite{hines}, but our proof is slightly different. Section 5 contains the main transcendence result and its proof (Theorem \ref{TeoGero01}). Section 6 points out another transcendence theorem. 

\begin{nota02}
For any (possibly finite) sequence in $\Za[i]$, $\sean$, we will write
\[
\left\langle a_0;a_1,a_2,\ldots\right\rangle:=a_0 + \cfrac{1}{a_1+\cfrac{1}{a_2+\cfrac{1}{\ddots}}}.
\]
In the following, only when $\sean$ is the Hurwitz continued fraction of some number we write $[a_0;a_1,a_2,\ldots]$ rather than $\left\langle a_0;a_1,a_2,\ldots\right\rangle$. By natural numbers, $\Na$, we mean the set of positive integers and we consider $\Na_0:=\Na\cup\{0\}$. By rational complex numbers we mean $\QU(i)$ and by irrational complex numbers $\Cx\setminus\QU(i)$. The real part of a complex number $z$ is $\Re(z)$ and its imaginary part is $\Im(z)$. For $A\subseteq \Cx$, the closure of $A$ is $\Cl A$ and the interior of $A$ is $A^{\circ}$ (both with respect to the usual topology), and for any $z\in\Cx$, $z+A :=\{z+a:a\in A \}$ and $zA :=\{za:a\in A \}$. Further notation is established along the text. 
\end{nota02}

\section{Hurwitz Continued Fractions}\label{SecHCF}

Denote by $[\cdot]:\Cx\to\Za[i]$ the function which assigns to each complex number its nearest Gaussian integers (choosing the one with largest real and imaginary part in case of ties). In symbols, if $\lfloor \cdot\rfloor:\RE\to\Za$ is the usual floor function, 
\[
\forall z\in\Cx \quad [z]:=\left\lfloor \Re(z) + \frac{1}{2}\right\rfloor + i\left\lfloor \Im(z) + \frac{1}{2}\right\rfloor;
\]
thus,
\[
\forall z\in\Cx \quad z - [z]\in \mfF:=\left\{w\in\Cx: -\frac{1}{2} \leq \Re (w), \Im(w) < \frac{1}{2}\right\}.
\]
By analogy with the Gauss map, we define on $\mfF^*:=\mfF\setminus \{0\}$ the function 
\[
T:\mfF^*\to\mfF, 
\qquad \forall z\in\mfF^* \quad
T(z)= \frac{1}{z} - \left[ \frac{1}{z}\right].
\]
Let $T^0:\mfF^*\to\mfF^*$ be the identitity map and $T^{n+1}:=T^n\circ T$ for $n\in\Na_0$. We associate to any $z\in\Cx$ a pair of possibly finite sequences, $(a_n)_{n\geq 0}$ and $(z_n)_{n\geq 0}$, via
\begin{align*}
z_0:=z, \quad & a_0:=[z_0], \nonumber\\
z_n := \frac{1}{T^{n-1}(z_0-a_0)},\quad & a_n:=[z_n], \nonumber
\end{align*}
for $n\in\Na$ as long as $T^{n-1}(z_0-a_0)\neq 0$.

\begin{figure}[ht]
\begin{center}
\includegraphics[scale=0.85,  trim={2.0in 6.75in 2.0in 0.75in},clip]{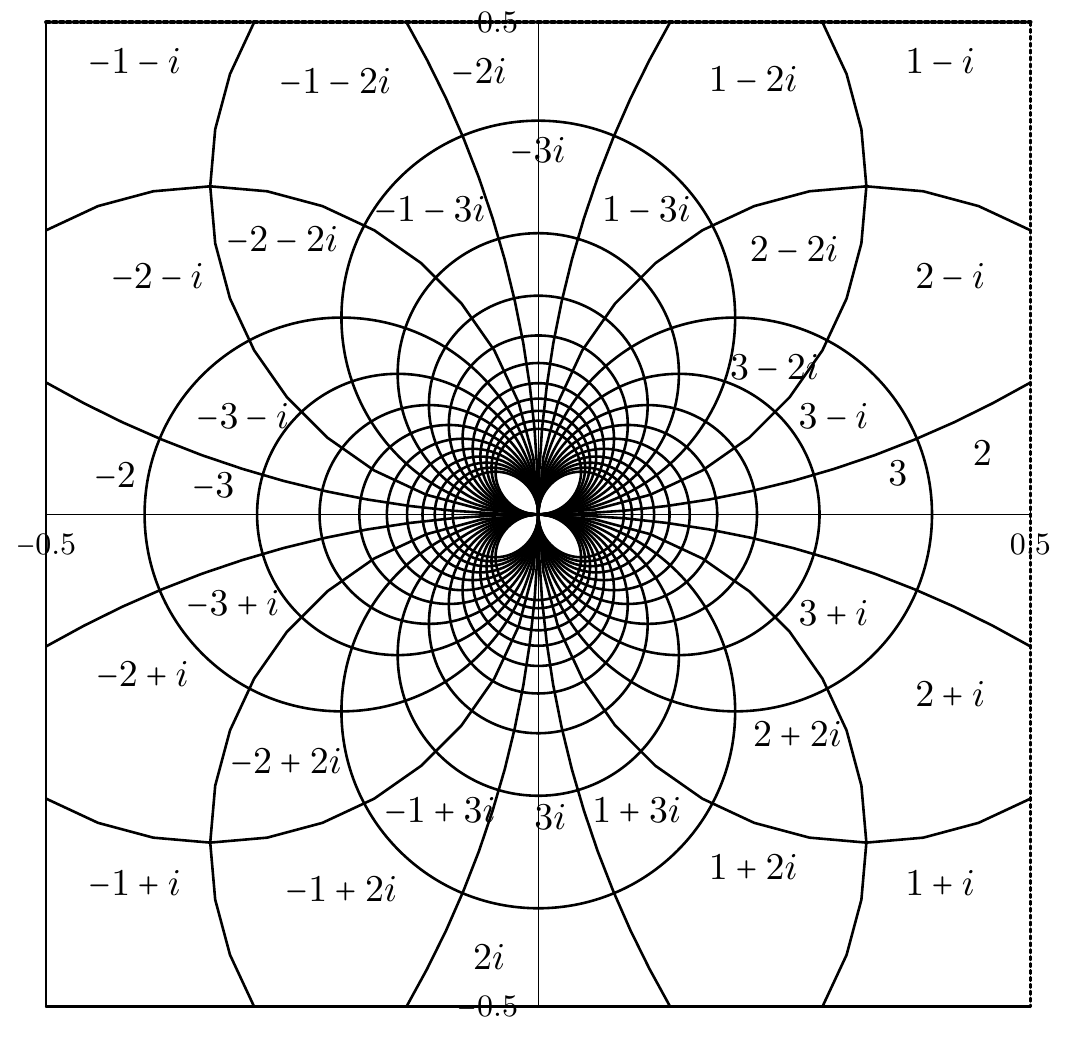}
\caption{Partition of $\mfF^*$ induced by $a_1$. \label{PartOfmfF}}
\end{center}
\end{figure}

Taking $a_0=0$ we can consider $a_1$ as a function from $\mfF^*$ to $\Za[i]$. The partition of $\mfF^*$ induced by the pre-images of $a_1$ is depicted in Figure \ref{PartOfmfF}. Note that
\begin{equation}\label{EcCpctDisj}
\mathfrak{K}:=\Cl\left\{ z\in\mfF: |a_1(z)|\geq \sqrt{8} \right\} \subseteq \mfF^{\circ}.
\end{equation}

\begin{def01}\label{DefHCF}
The \textbf{Hurwitz continued fraction} (HCF) of a complex number $z$ is the sequence $(a_n)_{n\geq 0}$ obtained by the above procedure. Let $\mathcal{I}\subseteq \Na_0$ be the set indexing $\sean$. As in \cite{daninog}, the $\mathcal{Q}$\textbf{-pair} of $z$ is the pair of sequences $\sepn$, $\seqn$ given by
\begin{equation*}
\begin{pmatrix}
p_{-2} & p_{-1} \\
q_{-2} & q_{-1}
\end{pmatrix}
=\begin{pmatrix}
0 & 1 \\
1 & 0
\end{pmatrix}, 
\quad
\forall n\in\mathcal{I} \quad
\begin{pmatrix}
p_{n} \\
q_{n}
\end{pmatrix}
=\begin{pmatrix}
p_{n-1} & p_{n-2} \\
q_{n-1} & q_{n-2} \\
\end{pmatrix} 
\begin{pmatrix}
a_{n} \\
1
\end{pmatrix}.
\end{equation*}
The terms of $( \tfrac{p_n}{q_n})_{n\geq 0}$ are the \textbf{HCF convergents} of $z$.
\end{def01}
Take $z\in\Cx$ and let $\sean$ be its HCF and $\sepn$, $\seqn$ its $\clQ$-pair. As expected, $\sean$ is infinite if and only if $z\in\Cx\setminus\QU(i)$ and in any case $z=[a_0;a_1,a_2,\ldots]$ (see \cite{daninog}, Theorem 3.7.). By the definition of $\sepn$ and $\seqn$, certain standard results are still true; for instance, for any valid $n\in\Na$ we have 
\begin{equation}\label{EcIdElem}
q_np_{n-1} - q_{n-1}p_n = (-1)^n, \quad
z= \frac{p_nz_{n+1}+p_{n-1}}{q_nz_{n+1}+q_{n-1}}
\end{equation}
(see \cite{daninog}, Proposition 3.3.).

Not every sequence in $\Za[i]$ is the HCF of a complex number. An infinite sequence of Gaussian integers is called \textbf{valid} if it is the HCF of some complex irrational number and we will denote the set of valid sequences by $\Omega^{\HCF}$. By a \textbf{valid prefix} we mean a finite sequence in $\Za[i]$ which is the prefix of a valid sequence. 

Some necessary conditions for a sequence to belong in $\Omega^{\HCF}$ follow immediately from the algorithm. 

\begin{figure}[ht] 
\begin{center}
\includegraphics[scale=1,  trim={2in 5in 2in 5in},clip]{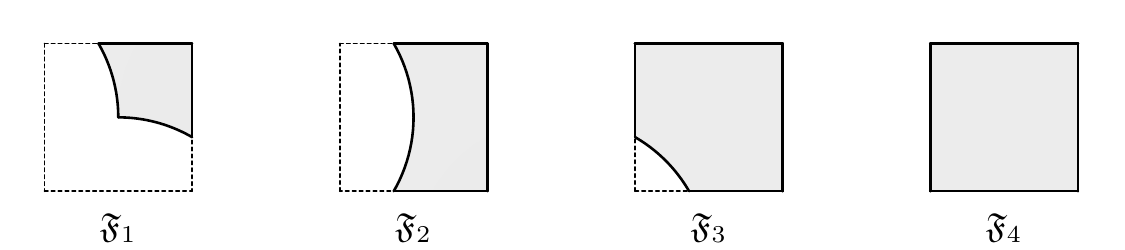}
\caption{\footnotesize{Feasible regions for $z_{n+1}^{-1}$.}\label{mfFCasos}}
\end{center}
\end{figure}

\begin{figure}[ht] 
\begin{center}
\includegraphics[scale=0.85,  trim={3.0in 8.8in 3.0in 0.7in},clip]{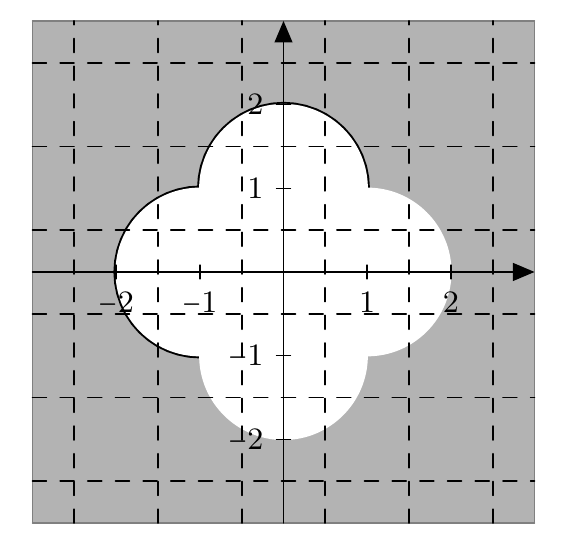}
\caption{ The set $\mfF^{-1}$.\label{FiguraF-1}}
\end{center}
\end{figure}

\begin{propo01}\label{PropoCondValida}
Let $(a_n)_{n\geq 0} \in\Omega^{\HCF}$ be the HCF of $z$, then for any admissible $n\in\Na$
\[
|z_n|\geq \sqrt{2}, \quad |a_n|\geq \sqrt{2}.
\]
We also have
\begin{equation}\label{CotaInfValidas}
\left\{ (b_n)_{n\geq 0}\in\Za[i]^{\Na_0}: \forall n\in\Na \quad |b_n| \geq \sqrt{8}\right\} \subseteq \Omega^{\HCF}.
\end{equation}
\end{propo01}
Take $z\in\mfF^*$, then $z_1\in \mfF^{-1}:=\{w^{-1}: w\in\mfF^*\}$ (see Figure \ref{FiguraF-1}). Suppose that $z_2$ exists. If $|a_1|\geq \sqrt{8}$, then the feasible maximal region for $z_2$ is again $\mfF^{-1}$. However, if $a_1=1+i$ ---for example--- then $z_2$ belongs to the set
\[
\left\{z\in\mfF^{-1}:   \Im (z) < \frac{1}{2},\; -\frac{1}{2}\leq \Re (z)  \right\}.
\]
Thus, determining whether a sequence is valid or not is more complicated than just checking a uniform lower bound. More generally, for a given $n\in\Na$ and a valid prefix $\bfa=(0;a_1,\ldots,a_n)$ define $\clC_n(\bfa):=\{z=[0;b_1,b_2,\ldots]\in\mfF^*: b_1=a_1,\ldots,b_n=a_n\}$. It can be shown that when $(T^n[\clC_n(\bfa)])^{\circ}\neq\vac$, there is some $j\in\Na$ such that $i^jT^n[\clC_n(\bfa)]\subseteq \mfF$ has the form of one of the sets depicted in Figure \ref{mfFCasos} (cfr. \cite{hiva18}, Section 3). We refer to the rules that determine whether a sequence is valid or not as the laws of succession.

The shift space associated to $(T,\mfF^*)$ cannot be modelled as a countable Markov shift. Indeed, assume there is an infinite the matrix $A$ characterizing\footnote{In the following sense: a sequence $(a_n)_{n\geq 1}$ in $\Za[i]$ belongs to $\Omega^{\HCF}$ if and only if $A_{a_na_{n+1}}=1$ for every $n\in\Na$.} $\Omega^{\HCF}$. On the one hand, the prefix $(0,1+2i,-2+2i,1+i)$ is not valid, so $A_{-2+2i,1+i}=0$. However, $(0,-2+2i,1+i)$ is a valid prefix, which would imply $A_{-2+2i,1+i}=1$, a contradiction. We can extend this observation to prefixes of arbitrary length. While for any $n\in\Na$ and any $\xi\in\Za[i]$ the sequences
\begin{align*}
&(\xi,\underbrace{-2+2i,2-2i,\ldots,2-2i,-2+2i}_{n \text{ repetitions of } -2+2i,2-2i},1+i), \nonumber\\
&(\xi,2-2i,\underbrace{-2+2i,2-2i,\ldots,2-2i,,-2+2i}_{n \text{ repetitions of }-2+2i,2-2i},1+i)\nonumber
\end{align*}
are valid prefixes, the sequences
\begin{align*}
&(\xi,1+2i,\underbrace{-2+2i,2-2i,\ldots,2-2i,-2+2i}_{n \text{ repetitions of }-2+2i,2-2i},1+i), \nonumber\\
&(\xi,-1+2i,2-2i,\underbrace{-2+2i,2-2i,\ldots,2-2i,,-2+2i}_{n \text{ repetitions of }-2+2i,2-2i},1+i)\nonumber
\end{align*}
are not. We can obtain even more examples using the symmetries of the HCF process. 

The lack of Markov structure is a significant difference between the Hurwitz continued fractions and its direct real analogue, the Ito-Tanaka $\tfrac{1}{2}$-continued fraction (cfr. \cite{tanito}). This feature complicates the study of finite sequences of the form $(a_j,a_{j-1},\ldots,a_1)$ where $(a_0,a_1,\ldots,a_{j-1},a_j)$ is a valid prefix. Such sequences appear naturally, because for any valid $\sean$ with $\clQ$-pair $\sepn$, $\seqn$ we have
\begin{equation}\label{EcMirror}
\forall n\in\Na \quad \frac{q_{n+1}}{q_n} = \left\langle a_{n+1};a_n,\ldots,a_1\right\rangle.
\end{equation}
Still, $(q_n)_{n\geq 0}$ has some desirable properties. For notational simplicity, we state the following lemma only for $\Cx\setminus\QU(i)$, although it also holds for $\QU(i)$.
\begin{lem01}\label{LemaDaniNogueira}
Let $z=[a_0;a_1,a_2,a_3,\ldots]$ be an element of $\Cx \setminus \QU(i)$ and let $\sepn, \seqn$ be its $\clQ$-pair. The following statements hold.
\begin{enumerate}[i.]
\item The sequence $(|q_n|)_{n\geq 0}$ is strictly increasing,
\item Set $\phi=\tfrac{1+\sqrt{5}}{2}$ , then
\[
\forall n\in\Na_0 \quad \frac{|q_{n+1}|}{|q_n|}> \phi \quad\text{or}\quad 
\frac{|q_{n+2}|}{|q_{n+1}|}> \phi,
\]
\item For every $n,k\in\Na$ we have
\[
|q_{n+k}|>\phi^{\left\lfloor  \frac{k}{2}\right\rfloor} |q_{n}|.
\]
\end{enumerate}
\begin{proof}
The first point is on page 195 of \cite{hur87}, the second is Corollary 5.3. of \cite{daninog} and the third follows directly from the previous two.
\end{proof}
\end{lem01}
In the following, we adopt the notation of \cite{lakein01}. For every $z\in\Cx$ and $\rho>0$,  $\Dx(z,\rho):=\{w\in\Cx: |w-z|<\rho\}$, $\overline{\Dx}(z,\rho):=\Cl\Dx(z,\rho)$, $\Ex(z,\rho):=\Cx\setminus \overline{\Dx}(z,\rho)$, and $C(z,\rho):=\{w\in\Cx:|w-z|=\rho\}$.
\section{Periodic Continued Fractions}

A. Hurwitz proved in \cite{hur87} an analogue to the Euler--Lagrange Theorem for his continued fractions (Theorem \ref{EuLaHCF} below). More than a century later, S.G. Dani and A. Nogueira generalized it to a larger family of complex continued fractions (\cite{daninog}, Corollary 4.3.).

\begin{teo01}[A. Hurwitz, \cite{hur87}]\label{EuLaHCF}
Let $z$ be an irrational complex number. The HCF of $z$ is ultimately periodic if and only if $z$ is quadratic over $\QU(i)$.
\end{teo01}

In \cite{galois}, \'Evariste Galois refined the Euler--Lagrange Theorem by showing that a quadratic irrational number $\alpha>1$ has a purely periodic regular continued fraction expansion if and only if its conjugate $\beta$ satisfies $-1<\beta<0$. Such numbers $\alpha$ are called \textbf{reduced}. We provide a similar result for HCF.

\begin{teo01}\label{PurelyPeriodicHCF}
Let $\xi=[a_0;a_1,a_2,\ldots]$ be a quadratic irrational over $\QU(i)$ and let $\eta \in\Cx$ be its conjugate over $\QU(i)$.
\begin{enumerate}[i.]
\item If $\sean$ is purely periodic, then $|\eta|< 1$.
\item If $|\xi|>1$, $\eta\in\mfF$ and $|a_n|\geq \sqrt{8}$ for every $n\in\Na_0$, then $\xi$ has purely periodic expansion. 
\item The conditions $\eta\in\mfF$ and $(\forall n\in\Na \; |a_n|\geq \sqrt{8})$ cannot be removed from the second point. In fact, there are infinitely many pairs $\xi,\eta$ such that
\begin{enumerate}[a.]
\item $\eta\not\in\mfF$, $|a_n|<\sqrt{8}$ for some $n$, and $\sean$ is not purely periodic,
\item $\eta\not\in\mfF$, $|a_n|\geq \sqrt{8}$ for all $n$, and $\sean$ is not purely periodic,
\item $\eta\in\mfF$, $|a_n|<\sqrt{8}$ for some $n$, and $\sean$ is not purely periodic.
\end{enumerate}
\end{enumerate}
\end{teo01}
If $\sean$ is a purely periodic sequence, we write $\sean=(\overline{a_0,\ldots,a_{m-1}})$ and we assume that $m\in\Na$ is minimal with respect to the property $a_{n}=a_{n+m}$ for every $n\in\Na_0$. A valid purely periodic sequence $(\overline{a_0,a_1,\ldots,a_{m-1}})$ is \textbf{reversible} if $(\overline{a_{m-1},\ldots,a_{1},a_0})$ is valid. Note that, by \eqref{CotaInfValidas}, every purely periodic sequence of Gaussian integers whose terms have absolute value at least $\sqrt{8}$ is reversible. 

Let us recall an elementary formula. If $\iota$ is the complex inversion: $\iota(z)=z^{-1}$ for $z\in\Cx\setminus\{0\}$; then, for all $z_0\in\Cx$ and all $\rho>0$ with $\rho\neq |z_0|$
\begin{equation}\label{EcInv}
\iota[C(z_0,\rho)] = C\left( \frac{\overline{z}_0}{|z_0|^2 - \rho^2}, \frac{\rho}{|\rho^2 - |z_0|^2|} \right).
\end{equation}
\begin{proof}[Proof of Theorem \ref{PurelyPeriodicHCF}]
\begin{enumerate}[i.]
\item Suppose that $\xi=[\overline{a_0;a_1,\ldots,a_{m-1}}]$. Let $j\in\Na$, then by \eqref{EcIdElem}
\[
\xi = \frac{p_{mj-1}\xi + p_{mj-2}}{q_{mj-1}\xi + q_{mj-2}},
\]
which implies
\[
q_{mj-1}\xi^2 + (q_{mj-2} - p_{mj-1})\xi - p_{mj-2}=0.
\]
Dividing by $q_{mj-1}$, we obtain a monic polynomial of second degree with coefficients in $\QU(i)$ which is satisfied by $\xi$:
\[
\xi^2 - \left(\frac{p_{mj-1}}{q_{mj-1}} -\frac{q_{mj-2}}{q_{mj-1}}  \right)\xi - \frac{p_{mj-2}}{q_{mj-1}}=0,
\]
so
\begin{equation}\label{EcEta+xi}
\eta + \xi = \frac{p_{mj-1}}{q_{mj-1}} -\frac{q_{mj-2}}{q_{mj-1}}.
\end{equation}
We conclude that
\begin{equation}\label{EcEta}
\lim_{j\to\infty} \frac{q_{mj-2}}{q_{mj-1}} = - \eta,
\quad\text{ hence }\quad
|\eta|\leq 1.
\end{equation}

In order to get $|\eta|<1$, we check two cases: $|a_{m-1}|\geq 2$ and $|a_{m-1}|= \sqrt{2}$. 

\paragraph{\S First case.} $|a_{m-1}|\geq 2$. Assume that $|q_{jm-1}/q_{jm-2}|<\phi$ holds for large $j$. For such $j$, Lemma \ref{LemaDaniNogueira} implies that $\phi\leq |q_{jm-2}/q_{jm-3}|$, so $|q_{jm-3}/q_{jm-2}|\leq \phi^{-1}$ and
\[
\left| \frac{q_{jm-1}}{q_{jm-2}} \right| 
=
\left| a_{m-1} + \frac{q_{jm-3}}{q_{jm-2}}  \right|
\geq 
2 - \frac{1}{\phi},
\]
hence $|\eta^{-1}|>1$. If $|q_{jm-1}/q_{jm-2}|\geq \phi$ holds for infinitely many $j$, then $|\eta^{-1}| >1$. 
\paragraph{\S Second case.}$|a_{m-1}|=\sqrt{2}$. By the symmetries of the HCF process, we do not lose generality by assuming that $a_{m-1}=1+i$. Then, $a_m=a_0$ must verify $\Im(a_0)\leq 0\leq \Re(a_0)$, so $m\geq 2$. Since $(|q_n|)_{n\geq 0}$ is strictly increasing, for all $j\in\Na_{\geq 2}$
\[
\frac{q_{jm-1}}{q_{jm-2}} = 1+i + \frac{q_{jm-3}}{q_{jm-2}}\in \Ex(0,1)\cap \Dx(1+i,1),
\]
which implies, by \eqref{EcInv}, that for all $j\in\Na_{\geq 4}$
\begin{equation}\label{Ec3bo}
\frac{q_{jm-2}}{q_{jm-3}} = a_{m-2} + \frac{q_{jm-4}}{q_{jm-3}}\in\Ex(-1+i,1)\cap \Ex(0,1) \cap \Dx(a_{m-2},1),
\end{equation}
that $\zeta:=\displaystyle \lim_{j\to\infty}\tfrac{q_{jm-3}}{q_{jm-2}}$ exists, and that $\zeta = -1-i-\eta^{-1}$. Assume for contradiction that $|\eta|=1$. Then, $\zeta\in C(-1-i,1)$ and
\begin{equation}\label{T1-Ec-Caso2-03}
\frac{1}{\zeta}\in C(-1+i,1) \cap \overline{\Dx}(a_{m-2},1),
\end{equation}
by \eqref{EcInv}. Thus, the only possibilities for $a_{m-2}$ are
\begin{equation}\label{LemaPPcaso201}
1+i, 2i, -1-i, -1+2i,  -1+3i, -2,  -2+i,  -2+2i, -3+i.
\end{equation}
(see Figure \ref{FiguraCasosFaciles}). The laws of succession exclude the options $1+i, \; 2i, -\; 1-i\;, -1+2i, \;-2, \; -2+i$. If we had $a_{m-2}=-1+3i$, then \eqref{T1-Ec-Caso2-03} would give
\[
\frac{1}{\zeta}\in C(-1+i,1)\cap \overline{\Dx}(-1+3i,1)=\{-1+2i\}
\]
and $\eta$ would belong to $\QU(i)$, a contradiction. $-3+i$ is discarded similarly. Then, we must have $a_{m-2}=-2+2i$. 

Since $(-2+2i, 1+i, -2+2i,1+i)$ is not valid, $m\geq 3$ holds. By \eqref{Ec3bo}, $\displaystyle\lim_{j\to\infty} \tfrac{q_{jm-4}}{q_{jm-3}}$ exists and belongs to $C(1-i,1)$. Arguing as above, we conclude that $a_{m-3}=2+2i$. We can repeat indefinitely the argument to obtain an infinite  sequence alternating between $-2+2i$ and $2+2i$, so $m\geq n$ for all $n\in\Na$, a contradiction. Therefore, $|\eta|<1$.

\begin{figure}[ht] 
\begin{center}
\includegraphics[scale=1,  trim={2.8in 4.25in 2.8in 4.05in},clip]{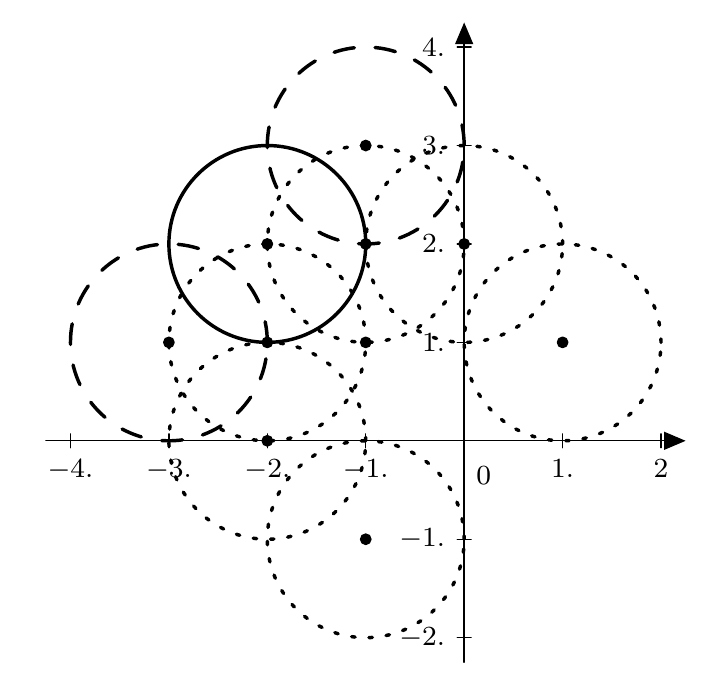}
\caption{The dotted circles are excluded by the laws of succession and the dashed circles by irrationality of $\zeta$. The solid circle is chosen. \label{FiguraCasosFaciles}}
\end{center}
\end{figure}
\item Let us keep the statement's notation. Conjugate over $\QU(i)$ the sequence $(\xi_j)_{j\geq 0}$ given by
\[
\xi_0:=\xi, \qquad \forall n\in\Na_0 \quad \xi_{n+1} = \frac{1}{\xi_n - a_n} \in\mfF
\]
and obtain
\[
\zeta_0:=\eta, \qquad \forall n\in\Na_0 \quad \zeta_{n+1} = \frac{1}{\zeta_n - a_n}.
\]
Let $(k_j)_{j\geq 0}$ be given by
\[
k_0=1, \qquad
\qquad 
\forall n\in\Na_0 \quad k_{n+1} = \sqrt{8} - \frac{1}{k_n}.
\]
It is not hard to show inductively that $k_n\to\sqrt{2}+1$ as $n\to\infty$, that $k_n>2$ for $n\geq 2$, and that $|\zeta_n|\leq k_n^{-1}$ for all $n\in\Na_0$.

For contradiction, assume that $\xi$ is not purely periodic. Write 
\[
\xi=[a_0;a_1,\ldots,a_n,\overline{a_{n+1},\ldots,a_{m+n}}]
\]
with $a_n\neq a_{n+m}$. Define $\xi'=[\overline{a_{n+1},\ldots,a_{m+n}}]=\xi_{n+1}$ and call $\eta'=\zeta_{n+1}$ its conjugate over $\QU(i)$. The proof of the first point tells us that
\[
-\frac{1}{\eta'} = [\overline{a_{m+n};a_{m+n-1},\ldots,a_{n+1}}].
\]

If we had $n\geq 2$, the left-most term in 
\[
\zeta_{n} - \zeta_{n+m} = a_{n} + \frac{1}{\zeta_{n+1}} - a_{n+m} -\frac{1}{\zeta_{n+m+1}} = a_{n} - a_{n+m}
\]
would belong to $\Dx(0,1)$ (since $|\zeta_n-\zeta_{n+m}|<k_n^{-1}+k_{n+m}^{-1}$), while the rightmost term would be a non-zero Gaussian integer. Therefore, either $n=0$ or $n=1$. Suppose that $n=1$. Conjugate $\xi_1=a_1+1/\xi'$ to get $\zeta_1=a_1+1/\eta'$ and
\begin{align}
-\zeta_1 &= -\frac{1}{\eta'} - a_1 \nonumber\\
&= [a_{m+1}-a_1;a_m,\ldots,a_2,\overline{a_{m+1},\ldots,a_2}] \in (a_{m+1}-a_1)+\mfF. \label{Ec-Zeta-1}
\end{align}

The inequality $|\zeta_1|\leq k_1^{-1}=(\sqrt{8}-1)^{-1}$, $a_{m+1}\neq a_1$, and \eqref{Ec-Zeta-1} give $|a_{m+1}-a_1|=1$. Assume that $a_{m+1}-a_1=1$ (the other cases are treated similarly). Since $\zeta_1\in\Dx(0,k_1)$,
\[
-\zeta_1-1\in\Dx(-1,(\sqrt{8}-1)^{-1})\cap \mfF; 
\]
hence, \eqref{EcInv} yields
\[
\frac{1}{-\zeta_1-1} \in \Dx\left( \frac{(\sqrt{8}-1)^2\sqrt{2}}{-8(\sqrt{2}-1)}, \frac{(\sqrt{8}-1)\sqrt{2}}{8(\sqrt{2}-1)}\right) \cap \mfF^{-1}.
\]
Call $c_0$ and $\rho_0$, respectively, the center and the radius of the last disc. Direct calculations give $0.7<\rho_0<0.8$ and $-1.5<c_0<-1.4$, so
\[
\left\{ a\in\Za[i] : \Dx(c_0;\rho_0) \cap \mfF^{-1} \cap (a+\mfF)\neq\vac\right\} = \{-2+i,-2,-2-i\}.
\]
Therefore, by \eqref{Ec-Zeta-1}, we would have $|a_m|\leq \sqrt{5}$, which contravenes $|a_m|\geq \sqrt{8}$. The only possibility left is $n=0$, but in this case we would have $\zeta_1=\eta'$ and
\begin{align*}
\eta&=a_0 + \frac{1}{\zeta_1} = a_0 - \frac{-1}{\eta'} \nonumber\\
&= a_0- a_{m} - [0;a_{m-1},\ldots,a_0,\overline{a_m,\ldots,a_0}] \in \mfF, \nonumber
\end{align*}
which implies $a_m=a_0$, a contradiction. Hence, $\sean$ is purely periodic. 

\item 
\begin{enumerate}[a.]
\item Take $M=M_1-iM_2, N=N_1-iN_2\in\Za[i]$ with $M_1,M_2,N_1,N_2\geq 2$. Since $(a_n)_{n\geq 1} = (\overline{M,1+i,N,2+4i})$ is valid and reversible, we can define the following pair of conjugate quadratic irrational numbers:
\[
\xi' := [\overline{M;1+i,N,2+4i}], 
\quad
-\frac{1}{\eta'} := [\overline{2+4i;N,1+i,M}].
\]
The conjugate of $\xi:=3+4i + 1/\xi'$ is
\begin{align*}
\eta &= 3+4i + \frac{1}{\eta'} = 3+4i - [\overline{2+4i;N,1+i,M}] \nonumber\\
&= [1;-N,-1-i,-M,\overline{-2-4i,-N,-1-i,-M}].
\end{align*}
Then, $\xi$ and $\eta$ have the next properties: $|\xi|>1$, $\xi$ has partial quotients with absolute value less than $\sqrt{8}$, $\eta\in\Dx(0,1)\setminus\mfF^*$ (because $\Re(-N)<0$), and $\xi$ does not have a purely periodic HCF. 

The core of the construction is to take a reversible sequence 
\[
(a_n)_{n\geq 1} = (\overline{a_1,\ldots,a_m}) \text{ with } \Re(a_{m-1})>0
\]
and pick $a_0=a_{m}+1$. 
\item The previous proof can be easily adapted.
\item The classical theory of regular continued fractions provides the first examples. Since reduced quadratic irrationals are dense in $(1,\infty)$, they are dense in $(1,2.5)$. Take a reduced quadratic irrational $\alpha$ with $1<\alpha<2.5$. The HCF of $\alpha=[a_0;a_1,a_2,\ldots]$ cannot be purely periodic, because $1<\alpha<1.5$ implies $a_0=1$ and $1.5<\alpha< 2.5$ implies $a_0=2, a_1<0$. However, the conjugate of $\alpha$ over $\QU$, and hence over $\QU(i)$, lies in $(-1,0)\subseteq \Dx$. An explicit example is the Golden Ratio $\phi$
\[
\phi=\frac{1+\sqrt{5}}{2} = [2;\overline{-3, 3}], \quad 
\frac{-1}{\phi}=\frac{1-\sqrt{5}}{2} = [-1; \overline{3, -3}].
\]
Although these examples do not show that $\sqrt{8}$ is best possible, they do provide a strategy to build infinitely many numbers showing it. Let $M\in\Za[i]\cap \mfF^{-1}$ satisfy $\Re(M)>0$. The sequence $(\overline{2+i;-2+i,M})$ is not valid, because any valid prefix of the form $(2+i;-2+i,M,2+i,-2+i,N)$ must have $\Re(N)<0$. However, $(\overline{M; -2+i,2+i})$ is valid, so, by \eqref{EcEta}, $-\left\langle 0;\overline{2+i,-2+i,M} \right\rangle$ converges to the conjugate of $[\overline{M;-2+i,2+i}]$. Then, we can define $\xi:=\xi(M):=\left\langle \overline{2+i;-2+i,M}\right\rangle$. Direct computations show that
\[
\xi = [2+i;-2+i,M+1,\overline{-2+i,2+i,M}].
\]
Yet, the conjugate of $\xi$, $\eta=-[0;\overline{M,-2+i,2+i}]$, is in $\mfF$. 
\end{enumerate}
\end{enumerate}
\end{proof}
\begin{rema1000}
The conjugate over $\QU(i)$ of $\xi=[\overline{a_0;a_1,\ldots,a_m}]$, $\eta$, satisfies 
\[
\eta=-\left\langle 0; \overline{a_m,\ldots,a_1,a_0}\right\rangle.
\]
This expansion is not necessarily a HCF. On certain occasions, we can compute the HCF of $\eta$ with  singularization\footnote{We borrow the term \textit{singularization} from \cite{ioskraa}.} identities such as 
\begin{align*}
a + \cfrac{1}{2+\cfrac{1}{b}} &= a+1 + \cfrac{1}{-2+\cfrac{1}{b+1}}, \nonumber\\
a+\cfrac{1}{1+i +\cfrac{1}{b}} &= a + (1-i) + \cfrac{1}{-(1+i) + \cfrac{1}{b+1-i}}  \nonumber
\end{align*}
for any $a,b\in\Cx$. A concrete example is 
\[
\xi=[\overline{5 + 6 i; -3 + 2 i, 2, 9 + 4 i}].
\]
The continued fraction of $\xi$ is not reversible. But we can get the HCF of $-\eta$ by reversing the period and applying the first singularization formula:
\[
\eta=-[0;\overline{10 + 4i, -2, -2 + 2i, 5 + 6 i}].
\]
\end{rema1000}

\section{Bounded Hurwitz Continued Fractions}

Minkowski's First Convex Body Theorem  (\cite{schm96}, Theorem 2B) yields a complex version of a famous corollary to Dirichlet's theorem on Diophantine approximation. Namely, a complex number $\zeta$ is irrational if and only if there are infinitely many co-prime Gaussian integers $p$ and $q$ such that
\begin{equation}\label{Ec4pi}
\left| \zeta - \frac{p}{q}\right| < \frac{4}{\pi}\frac{1}{|q|^2}.
\end{equation}
A complex number is badly approximable if \eqref{Ec4pi} cannot be substantially improved.

\begin{def01}\label{DefBAD}
A complex number $\zeta$ is \textbf{badly approximable} if there exists a constant $C>0$ such that for every $p,q\in\Za[i]$ with $q\neq 0$
\[
\left| \zeta - \frac{p}{q}\right| > \frac{C}{|q|^2}.
\]
$\bad_{\Cx}$ denotes the set of badly approximable complex numbers.
\end{def01}

$\bad_{\Cx}$ shares some properties with its real counterpart. For instance, it has Lebesgue measure $0$, it is $\frac{1}{2}$-winning in the sense of Schmidt games (\cite{dodkrist}, Theorem 5.2.), and hence, it has full Hausdorff dimension. We can also characterize $\bad_{\Cx}$ in terms of Hurwitz continued fractions.

\begin{teo01}\label{BadCCharac}
The following equality holds
\[
\bad_{\Cx} = \left\{z=[a_0;a_1,a_2,\ldots] \in\Cx\setminus\QU(i): \exists M>0 \quad \forall n\in\Na_0 \quad |a_n|\leq M\right\}.
\]
\end{teo01}

Most of the standard argument (\cite{khin}, Theorem 23) used in the real version of Theorem \ref{BadCCharac} also works in our context. However, we must establish some approximation properties of HCF. Take any $\zeta\in\Cx$ and $p,q\in\Za[i]$ co-prime with $q\neq 0$. We say that $p/q\in\QU(i)$ is a \textbf{good approximation} to $\zeta$ if
\[
|q\zeta - p| = \min\left\{| q'\zeta-p'|: q',p'\in\Za[i] \quad |q'| \leq |q|\right\}.
\]
We say that $p/q\in\QU(i)$ is a \textbf{best approximation} to $\zeta$ if
\[
\forall p',q'\in \Za[i] \qquad
|q'|<|q| \quad\implies\quad
|q\zeta-p|< |q'\zeta-p'|.
\] 
The next result is Theorem 1 of \cite{lakein01}.

\begin{teo01}[R. Lakein, \cite{lakein01}]\label{Lakein72Teo02}
For all $\zeta\in\Cx$ every HCF convergent of $\zeta$ is a good approximation to $\zeta$. Moreover, for almost every $\zeta\in\Cx$ (with respect to the Lebesgue measure) every HCF convergent of $\zeta$ is a best approximation to $\zeta$.
\end{teo01}

\begin{lem01}\label{PropoLemaBAD}
Let $\zeta=[a_0;a_1,a_2,\ldots]$ be in $\Cx\setminus\QU(i)$ and let $\sepn$, $\seqn$ be its $\clQ$-pair; then, writing $\gamma=\left( 1-\tfrac{\sqrt{2}}{2} \right)^{-1}$,
\begin{equation}\label{Ec-LemaBAD}
\forall n\in\Na \quad
\frac{1}{\left(|\zeta_{n+1}| + 1\right)|q_n|^2}<\left| \zeta - \frac{p_n}{q_n}\right| \leq \frac{\gamma}{|q_nq_{n+1}|}.
\end{equation}
\end{lem01}
\begin{proof}
Both of the equations in \eqref{EcIdElem} give
\[
\zeta = \frac{p_n\zeta_{n+1}+p_{n-1}}{q_n\zeta_{n+1}+q_{n-1}} = \frac{p_n}{q_n} + \frac{(-1)^n}{q_n^2\left( \zeta_{n+1} + \frac{q_{n-1}}{q_n}\right)}. 
\]
The lower bound in \eqref{Ec-LemaBAD} follows immediately from the strict monotonicity of $(|q_n|)_{n\geq 0}$. For the upper bound, we use $\zeta_{n+1}=a_{n+1}+\zeta_{n+2}^{-1}$ to obtain
\begin{equation}\label{Ec-LemaBadLD}
\zeta - \frac{p_n}{q_n} 
= \frac{(-1)^n}{q_n(q_n\zeta_{n+1}+q_{n-1})} 
= \frac{(-1)^n}{q_nq_{n+1}\left(1 + \frac{q_n}{q_{n+1}}\frac{1}{\zeta_{n+2}}\right)}.
\end{equation}

Since $(|q_n|)_{n\geq 0}$ is strictly increasing and $\zeta_{n+2}^{-1}\in\mfF\subseteq \overline{\Dx}(0,2^{-\frac{1}{2}})$, 
\[
\left| 1+ \frac{q_n}{q_{n+1}}\frac{1}{\zeta_{n+2}}\right| \geq 1 - \frac{\sqrt{2}}{2}. 
\]
The previous inequality and \eqref{Ec-LemaBadLD} yield the result. 
\end{proof}

\begin{proof}[Proof of Theorem \ref{BadCCharac}]
We start with $\supseteq$. Let $\zeta=[a_0;a_1,a_2,\ldots]$ be a complex number whose partial quotients are bounded by $M>0$. Then,
\[
\forall n\in\Na \quad \left| \zeta_n\right| = \left| a_n+\frac{1}{\zeta_{n+1}} \right| \leq M+1
\]
and, by \eqref{Ec-LemaBAD}, for $c_1=(M+2)^{-1}$ we obtain
\[
\forall n\in\Na \qquad  \frac{c_1}{|q_n|^2} \leq \left| \zeta - \frac{p_n}{q_n}\right|.
\]
Now, take $p/q\in\QU(i)$ in its lowest terms and $n\in\Na$ such that $|q_{n-1}| < |q|\leq |q_n|$. By Theorem \ref{Lakein72Teo02}, $|q_n\zeta - p_n|\leq |q\zeta - p|$, so 
\[
\left| \zeta - \frac{p}{q}\right| \geq \frac{|q_n|}{|q|} \,\left| \zeta - \frac{p_n}{q_n} \right| 
>\frac{c_1}{|q|^2} \frac{|q_{n-1}|^2}{|q_n|^2} 
= \frac{c_1}{|q|^2} \frac{1}{\left|a_n + \frac{q_{n-2}}{q_{n-1}}\right|^2} 
\geq \frac{c_1}{|q|^2(M+1)^2}.
\]
Setting $C=(M+1)^{-1}(M+2)^{-2}$ we conclude that $\zeta\in\bad_{\Cx}$.

In order to show $\subseteq$, take $\zeta\in\bad_{\Cx}$ and let $C=C(\zeta)>0$ be the constant from Definition \ref{DefBAD}. Set $\gamma$ as in Lemma \ref{PropoLemaBAD} and $M=\tfrac{\gamma}{C}+1$, then for any $n\in\Na$
\begin{align*}
\frac{C}{|q_n|^2} \leq \frac{\gamma}{|q_nq_{n+1}|} \quad &\implies\quad \left| a_{n+1} + \frac{q_{n-1}}{q_n}\right| = \left| \frac{q_{n+1}}{q_n}\right| \leq \frac{\gamma}{C}\nonumber\\
&\implies\quad |a_{n+1}|\leq M \nonumber
\end{align*}
\end{proof}

We can restate Theorem \ref{TeoBoGru} as follows.
\begin{teo01}
There are badly approximable algebraic complex numbers of arbitrary even degree over $\QU(i)$.
\end{teo01}
Recently, Robert Hines gave another proof of Theorem \ref{BadCCharac} (\cite{hines}, Theorem 1). While his argument relies on R. Lakein's work too, it avoids \eqref{Ec-LemaBAD}.
\section{Transcendental Complex Numbers}

Recall that the \textbf{length} of a finite word $\bfx$, $|\bfx|$, is the number of terms it comprises. Let $\clA\neq \vac$ be a finite alphabet and let $\bfa$ be an infinite word over $\clA$. As noted in \cite{bugkim}, $\rep(\bfa)<+\infty$ is equivalent to the existence of three sequences of finite words in $\clA$, $(W_n)_{n\geq 1}$, $(U_n)_{n\geq 1}$, and $(V_n)_{n\geq 1}$, such that
\begin{enumerate}[i.]
\item For every $n$ the word $W_nU_nV_nU_n$ is a prefix of $\bfa$,
\item The sequence $((|W_n|+|V_n|)/|U_n|)_{n\geq 1}$ is bounded above,
\item The sequence $(|U_n|)_{n\geq 1}$ is strictly increasing.
\end{enumerate}

Our main result is the following theorem.
\begin{teo01}\label{TeoGero01}
Let $\bfa=(a_j)_{j\geq 0}$ be a non-periodic, valid sequence such that
\[
\rep \bfa < + \infty.
\]
Define $\zeta=[a_0;a_1,a_2,a_3,\ldots]$ and let $(W_n)_{n\geq 0},(U_n)_{n\geq 0},(V_n)_{n\geq 0}$ be as above. 
\begin{enumerate}[i.]
\item If $\displaystyle\liminf_{n\to\infty} |W_n| < + \infty$, then $\zeta$ is transcendental.
\item If $\displaystyle\liminf_{n\to\infty} |W_n| = + \infty$ and $|a_n|\geq \sqrt{8}$ for every $n\in\Na$, then $\zeta$ is transcendental.
\end{enumerate}
\end{teo01}

\subsection{Preliminary Results}

We will keep the notation of Theorem \ref{TeoGero01} until the end of Section 5. Let $\sepn$, $\seqn$ be the $\mathcal{Q}$-pair of $\zeta$. Define the sequences of non-negative integers $(w_n)_{n\geq 1}$, $(u_n)_{n\geq 1}$, $(v_n)_{n\geq 1}$ by
\[
\forall n\in\Na \quad w_n := |W_n|, \quad u_n := |U_n|, \quad v_n := |V_n|.
\]

\begin{lem01}\label{LemmaTeoGero}
Let $\psi:=\left( \frac{1+\sqrt{5}}{2}\right)^{\frac{1}{2}}$. There exists some $\veps>0$ such that
\[
\forall n\in\Na \quad \psi^{u_n}\geq |q_{w_n}q_{w_n+u_n+v_n}|^{\veps}.
\]
\end{lem01}
\begin{proof}
In view of the the boundedness of $\sean $ and $((v_n+w_n)/u_n)_{n\geq 0}$, we can define the real numbers
\[
M := 1+\sup_{n\in\Na} |q_n|^{\frac{1}{n}}, \quad 
N := 2 + \sup_{n\in\Na} \frac{2w_n + v_n}{u_n}.
\]
Since $\psi>1$, for every $n\in\Na$
\[
\psi = (M^{w_n}M^{w_n+u_n+v_n})^{\frac{\log\psi}{(2w_n+u_n+v_n)\log M}} > |q_{w_n}q_{u_n+v_n+w_n}|^{\frac{\log\psi}{(2w_n+u_n+v_n)\log M}},
\]
and $\veps=\frac{\log\psi}{N\log M}$ works:
\[
\psi^{u_n} > \left(|q_{w_n}q_{u_n+v_n+w_n}|^{\frac{u_n}{2w_n+u_n+v_n}}\right)^{\frac{\log\psi}{\log M}}
> |q_{w_n}q_{u_n+v_n+w_n}|^{\frac{\log\psi}{N\log M}}.
\]
\end{proof}

As in \cite{bug13autom}, our main tool is an adequate version of Schmidt's Subspace Theorem. We will use a particular case of the corresponding result for number fields (cfr. \cite{schm76}, Theorem 3). Before stating it, we require some notation. If $z\in\Cx$, then $\overline{z}$ is its complex conjugate. If $\bfz=(z_1,\ldots,z_k)\in\Cx^k$, then $\overline{\bfz}:=(\overline{z_1}, \ldots,\overline{z_k})\in\Cx^k$ and $\| \bfz\|_{\infty}:=\max\{|z_1|,\ldots,|z_k|\}$. We denote the zero vector by $\bfzo$.

\begin{teo01}[Schmidt's Subspace Theorem]\label{ScSsThQi}
Let $\scL_1$,$\ldots$,$\scL_m$, $\scM_1$,$\ldots$,$\scM_m$ be two sets of $m\in\Na$ linearly independent linear forms in $m$ variables. Suppose that all the forms have algebraic coefficients. Then, for any $\veps>0$ there are finitely many proper subspaces of $\QU(i)^m$, $T_1,\ldots,T_k$, such that for all $\bfbe \in \Za[i]^m$, $\bfbe\neq\bfzo$,
\[
\left| \prod_{j=1}^m \scL_j(\bfbe)\right|\left| \prod_{j=1}^m \scM_j(\overline{\bfbe})\right| \leq \frac{1}{\|\bfbe\|_{\infty}^{\veps}} 
\quad\implies\quad
\bfbe\in\bigcup_{j=1}^k T_j.
\]
\end{teo01}
\subsection{Proof of Theorem \ref{TeoGero01}}

Let us restate the two cases of Theorem \ref{TeoGero01} in a simpler way. In the first one, rather than $\displaystyle\liminf_n w_n<+\infty$, we may assume that $(w_n)_{n\geq 1}$ is constant after taking an appropriate sub-sequence. In the second one, we replace $\liminf_n w_n=+\infty$ by two conditions: $(w_n)_{n\geq 1}$ is strictly increasing and $a_{w_n}\neq a_{w_n+u_n+v_n}$ for all $n\in\Na$ (see \cite{bug13autom}, p.1013).

There is no loss of generality if we assume $a_0=0$. Indeed, in general, if $\bfx$ and $\bfy$ are two infinite words over a finite alphabet differing only in the first term, then $r(n-1,\bfx)\leq r(n,\bfy)\leq r(n+1,\bfx)$ for all $n\in\Na$ (see Definition \ref{DefRepExp}). Therefore, $\rep(\bfx)=\rep(\bfy)$. Also, $\zeta$ is transcendental if and only if $\zeta-k$ is transcendental for all $k\in\Za[i]$. 

Note that $a_0=0$ is equivalent to $\zeta\in\mfF$, which implies
\[
\forall n\in\Na_0 \quad |p_n|< |q_n|.
\]
We will use the previous inequality and the first part of Lemma \ref{LemaDaniNogueira} without reference.

Assume, for contradiction, that $\zeta$ is algebraic. Observe that, since $\bfa$ is not periodic, $3\leq [\QU(\zeta,i):\QU(i)]<+\infty$. 

\paragraph{\S.$(w_n)_{n\geq 1}$ is constant.} We can suppose that $w_n=0$ for all $n$, because the transcendence of any member of $\{[a_n;a_{n+1},a_{n+2},\ldots]:n\in\Na\}$ implies the transcendence of the rest.

Write $s_n=u_n+v_n$ for every $n\in\Na$. 
\subparagraph{\textsc{i.}} Define the sequence of irrational quadratic numbers $\left(\zeta^{(n)}\right)_{n\geq 1}$ by
\begin{equation}\label{EcCaso01Def-01}
\forall n\in\Na \qquad \zeta^{(n)} = [0;b_{n,1},b_{n,2},b_{n,3},\ldots]
 = [0;U_nV_nU_nV_nU_nV_n\ldots].
\end{equation}
The second equality should be understood as follows. For a given $n\in\Na$ and any $j\in\Na$, let $r\in\Na$ be such that
\[
1\leq r\leq s_n, \qquad r\equiv j \pmod{s_n}.
\]
Then, writing $U_nV_n=d_1d_2\ldots d_{s_n}$, we set $b_{n,j}=d_r$.

We might need to take a sub-sequence of $(U_n)_{n\geq 1}, (V_n)_{n\geq 1}$ to guarantee $0\overline{U_nV_n}\in\Omega^{\HCF}$ for each $n$. We sketch an explanation and leave the details to the reader. In general, let $X=x_1\ldots x_m$ and $Y=y_1\ldots y_k$ be two non-empty finite words on $\Za[i]$ such that $0XYX$ is a valid prefix but $0XYXY$ is not. Since only the the first appearance of $x_my_1$ is allowed, we must have $|x_m|\in\{\sqrt{5},\sqrt{8}\}$. After studying separately each case, we conclude that for some fixed $l\in\{1,2,3,4\}$ the other terms of $X$ alternate between $i^l(2+2i)$ and $i^l(-2+2i)$ and that for some $j$ we have $|y_j|=\sqrt{5}$. If there were infinitely many numbers $n$ for which $0\overline{U_nV_n}$ is not valid, we would obtain a contradiction with the aid of the the three conditions satisfied by $(U_n)_{n\geq 1}$ and $(V_n)_{n\geq 1}$. Thus, $(\zeta^{(n)})_{n\geq 1}$ is well defined.

Let $m$ be a natural number. Since the first $u_m+s_m$ terms of the HCF of $\zeta$ and $\zeta^{(m)}$ coincide, Lemma \ref{PropoLemaBAD} implies that for some absolute constant $\kappa_0>0$ we have
\[
\left| \zeta - \zeta^{(m)}\right| < \frac{\kappa_0}{|q_{u_m+s_m}|^2}.
\]
As in the proof of Theorem \ref{PurelyPeriodicHCF}, $\zeta^{(m)}$ satisfies the polynomial
\[
P_m(X)  := q_{s_m-1}X^2 + \left(q_{s_m}-p_{s_m-1}\right)X - p_{s_m}.
\]
Therefore, there are constants $\kappa_1=\kappa_1(\zeta)>0$, $\kappa_2=\kappa_2(\zeta)>0$ such that
\begin{align}
\forall n\in\Na\quad 
| P_n(\zeta)|  &= \left| P_n(\zeta)  - P_n(\zeta^{(n)}) \right| \nonumber\\
&= \left| q_{s_n-1}\left(\zeta - \zeta^{(n)}\right)\left(\zeta + \zeta^{(n)}\right) + \left(q_{s_n}-p_{s_n-1}\right)\left(\zeta - \zeta^{(n)}\right) \right| \nonumber\\
&\leq \kappa_1\left( |q_{s_n}|\left| \zeta - \zeta^{(n)}\right| + 2|q_{s_n}|\left| \zeta - \zeta^{(n)}\right| \right) \nonumber\\
&\leq \kappa_2 \frac{|q_{s_n}|}{|q_{u_n+s_n}|^2} \label{EcCaso1-01}.
\end{align}
\subparagraph{\textsc{ii.}}  Let $(\bfx_n)_{n\geq 1}$ be the sequence in $\Za[i]^4$ given by
\[
\forall n\in\Na \quad 
\bfx_n = (-q_{s_n-1},-p_{s_n-1}, q_{s_n}, p_{s_n}), \quad\text{so}\quad  \| \bfx_n \|_{\infty} = |q_{s_n}|.
\]
Define two sets of independent linear forms $\scL_1^1$, $\scL_2^1$, $\scL_3^1$, $\scL_4^1$, $\scM_1^1$, $\scM_2^1$, $\scM_3^1$, $\scM_4^1$ in the variables $\bfX=(X_1,X_2,X_3,X_4)$ and $\widetilde{\bfX}=(\widetilde{X}_1,\widetilde{X}_2,\widetilde{X}_3,\widetilde{X}_4)$, respectively, by
\begin{align*}
\scL_1^1(\bfX) &= \zeta^2 X_1 - \zeta(X_2 + X_3) + X_4,  &&\scM_1^1(\widetilde{\bfX}) = \overline{\zeta}^2 \widetilde{X}_1 - \overline{\zeta}(\widetilde{X}_2 + \widetilde{X}_3) + \widetilde{X}_4,\nonumber\\
\scL_2^1(\bfX) &= \zeta X_1 - X_2,  &&\scM_2^1(\widetilde{\bfX})= \overline{\zeta} \widetilde{X}_1 - \widetilde{X}_2,\nonumber\\
\scL_3^1(\bfX) &= \zeta X_1 - X_3,  &&\scM_3^1(\widetilde{\bfX})= \overline{\zeta} \widetilde{X}_1 - \widetilde{X}_3, \nonumber\\
\scL_4^1(\bfX) &= X_1, &&\scM_4^1(\widetilde{\bfX})= \widetilde{X}_1.
\end{align*}

Set $\psi=\left( \frac{1+\sqrt{5}}{2}\right)^{\frac{1}{2}}$. By Lemma \ref{PropoLemaBAD}, \eqref{EcCaso1-01}, and the third point of Lemma \ref{LemaDaniNogueira}, there is some $\kappa_3=\kappa_3(\zeta)>0$ such that for any $n\in\Na$
\begin{align*}
\left| \scL_1^1 \scL_2^1 \scL_3^1 \scL_4^1(\bfx_n)\right| &=|P_n(\zeta)||\zeta q_{s_n-1} - p_{s_n-1}||\zeta q_{s_n-1} + q_{s_n}||q_{s_n-1}| \nonumber\\
& \leq \kappa_3 \frac{|q_{s_n}|^2}{|q_{s_n+u_n}|^2} <\kappa_3\frac{\psi^2}{\psi^{2u_n}}. \nonumber
\end{align*}

So there are constants $\kappa_4=\kappa_4(\zeta)>0$ and $\veps>0$ (by Lemma \ref{LemmaTeoGero}) such that
\begin{equation*}
\forall n\in\Na \quad \left| \scL_1^1 \scL_2^1 \scL_3^1 \scL_4^1 (\bfx_n)\right| \leq \frac{\sqrt{\kappa_4}}{|q_{s_n}|^{\veps}} = 
\frac{\sqrt{\kappa_4}}{\|\bfx_n\|_{\infty}^{\veps}}.
\end{equation*}
For any $j\in\{1,2,3,4\}$ and any $n\in\Na$ we have $|\scM^1_j(\overline{\bfx}_n)|=|\scL^1_j(\bfx_n)|$, then
\[
\forall n\in\Na \quad  \left| \scL_1^1 \scL_2^1 \scL_3^1 \scL_4^1 (\bfx_n)\right|\left| \scM_1^1 \scM_2^1 \scM_3^1 \scM_4^1 (\overline{\bfx_n})\right| \leq \frac{\kappa_4}{\|\bfx_n\|_{\infty}^{2\veps}}.
\]
Since $\|\bfx_n\|\to\infty$ as $n\to\infty$, for large $n$ we have $\kappa_4/\|x_n\|^{\veps}<1$ and
\begin{equation}\label{PnAlfaAux03}
\left| \scL_1^1 \scL_2^1 \scL_3^1 \scL_4^1 (\bfx_n)\right|\left| \scM_1^1 \scM_2^1 \scM_3^1 \scM_4^1 (\overline{\bfx_n})\right| \leq \frac{1}{\|\bfx_n\|_{\infty}^{\veps}}.
\end{equation}
Then, by Theorem \ref{ScSsThQi}, there exist a vector $\bfx=(x_1,x_2,x_3,x_4)\in\Za[i]^4$, $\bfx\neq \bfzo$, and an infinite set $\clN_1\subseteq\Na$ satisfying
\[
\forall n\in\clN_1 \quad x_1q_{s_n-1}+x_2p_{s_n-1}+x_3q_{s_n}+x_4p_{s_n}=0
\]
Dividing the last expression by $q_{s_n-1}$ we get
\begin{equation}\label{Ec22Tilde}
\forall n\in\clN_1 \quad x_1+x_2\frac{p_{s_n-1}}{q_{s_n-1}}+x_3\frac{q_{s_n}}{q_{s_n-1}}+x_4\frac{p_{s_n}}{q_{s_n}}\frac{q_{s_n}}{q_{s_n-1}}=0.
\end{equation}
Since $\zeta\in\Cx\setminus\QU(i)$, $x_3,x_4\in\Za[i]$ cannot be both $0$; thus, we can define
\[
\xi = \lim_{\substack{n\to\infty\\ n\in\clN_1}} \frac{q_{s_n}}{q_{s_n-1}} = -\frac{\zeta x_2 + x_1}{\zeta x_4 + x_3}.
\]
Let us show that $\xi\in\Cx\setminus\QU(i)$. By \eqref{Ec-LemaBAD}, for some $C_1=C_1(\bfx,\zeta)>0$, $C_2=C_2(\bfx,\zeta)>0$, every $n\in \clN_1$ verifies
\begin{align}
\left| \xi - \frac{q_{s_n-1}}{q_{s_n}}\right| &= \left| \frac{x_1+\zeta x_3}{x_2 + \zeta x_4} - \frac{x_1+\frac{p_{s_n}}{q_{s_n}} x_3}{x_2 + \frac{p_{s_n-1}}{q_{s_n-1}} x_4}  \right| \nonumber\\
&\leq \left| \frac{x_1}{x_2+\zeta x_4} - \frac{x_1}{x_2 + x_4\frac{p_{s_n-1}}{q_{s_n-1}}}\right| 
+ 
|x_3| \left| \frac{\zeta}{x_2+x_4\zeta} - \frac{\frac{p_{s_n}}{q_{s-n}}}{x_2 + \frac{p_{s_n-1}}{q_{s_n-1}}x_4}\right|
\nonumber\\
&\leq C_1\left( \left| \zeta - \frac{p_{s_n}}{q_{s_n}}\right| + \left| \zeta - \frac{p_{s_n-1}}{q_{s_n-1}}\right| + \left| \frac{p_{s_n-1}}{q_{s_n-1}} - \frac{p_{s_n}}{q_{s_n}} \right|\right)\nonumber\\
&\leq \frac{C_2}{|q_{s_n}q_{s_n-1}|}. \label{approx-xi}
\end{align}

Assume that $\xi=a/b$ with $a,b\in\Za[i]$ co-prime. Since $q_{s_j-1}$ and $q_{s_j}$ are co-prime for all $j\in\Na$, for large $n\in\clN_1$ we would have
\[
\frac{1}{|bq_{s_{n-1}}|} \leq \left| \xi - \frac{q_{s_n-1}}{q_{s_n}}\right| \leq \frac{C_2}{|q_{s_n-1}q_{s_n}|},
\]
which implies $|q_{s_n-1}|\leq C_2|b|$. However, this contradicts the unboundedness of $(|q_{s_j-1}|)_{j\in\clN_1}$; hence, $\xi\in\QU(i,\zeta)$ is irrational. Furthermore, letting $n\to\infty$ along $\clN_1$ in \eqref{Ec22Tilde} and recalling that $[\QU(i,\zeta):\QU(i)]\geq 3$, we get $\xi\neq \zeta$.

\subparagraph{\textsc{iii.}} Consider the linear forms $\scL_1^2$, $\scL_2^2$, $\scL_3^2$, $\scM_1^2$, $\scM_2^2$, $\scM_3^2$ in the variables $\bfY=(Y_1,Y_2,Y_3)$ and $\widetilde{\bfY}=(\widetilde{Y}_1,\widetilde{Y}_2,\widetilde{Y}_3)$, respectively, given by
\begin{align*}
\scL_1^2(\bfY) &= \xi Y_1 - Y_2, && \scM_1^2(\widetilde{\bfY}) = \overline{\xi} \widetilde{Y}_1 - \widetilde{Y}_2, \nonumber\\
\scL_2^2(\bfY) &= \zeta Y_1 - Y_3, && \scM_2^2(\widetilde{\bfY}) = \overline{\zeta} \widetilde{Y}_1 - \widetilde{Y}_3, \nonumber\\
\scL_3^2(\bfY) &= Y_2, && \scM_3^2(\widetilde{\bfY}) = \widetilde{Y}_2 .\nonumber
\end{align*}

Define the sequence $(\bfy_n)_{n\geq 0}$ in $\Za[i]^3$ by
\[
\forall n\in\Na \quad
\bfy_n:= (q_{s_n},q_{s_n-1},p_{s_n}), \quad\text{ so}\quad \|\bfy_n\|_{\infty}= |q_{s_n}|.
\]
By Lemma \ref{PropoLemaBAD} and \eqref{approx-xi}, for some $\kappa_5=\kappa_5(\zeta)>0$ and every $n\in\clN_1$
\begin{align*}
\left|\scM_1^2\scM_2^2\scM_3^2(\overline{\bfy_n}) \right| &=\left|\scL_1^2\scL_2^2\scL_3^2(\bfy_n) \right| \nonumber\\
&= |(q_{s_n}\xi - q_{s_n-1})(q_{s_n}\zeta - p_{s_n})q_{s_n-1}| \nonumber\\
& \leq \kappa_5 \frac{|q_{s_n-1}|}{|q_{s_n-1}q_{s_n}|} = 
 \frac{\kappa_5}{\| \bfy_n \|_{\infty}}. \nonumber
\end{align*}
Thus, arguing as above, for every large $n\in\clN_1$
\[
\left|\scL_1^2\scL_2^2\scL_3^2 (\bfy_n) \right| \left|\scM_1^2\scM_2^2\scM_3^2 (\overline{\bfy_n} )\right| \leq \frac{1}{\|\bfy_n\|_{\infty}}.
\]
Theorem \ref{ScSsThQi} yields the existence of some $\bfzo\neq\bfy=(y_1,y_2,y_3)\in\Za[i]^3$ and of an infinite set $\clN_2\subseteq\clN_1$ such that
\[
\forall n\in\clN_2 \quad q_{s_n}y_1 + q_{s_n-1}y_2 + p_{s_n}y_3 = 0.
\]
Dividing by $q_{s_n}$ and taking the limit when $n\to\infty$ along $\clN_2$ we obtain
\begin{equation}\label{EcAlfaBeta01}
y_1 + \xi y_2 + \zeta y_3 =0 \quad\text{ and }\quad y_2y_3\neq 0,
\end{equation}
because $\xi,\zeta\in\Cx\setminus\QU(i)$.
\subparagraph{\textsc{iv.}} Define $\scL_1^3$, $\scL_2^3$, $\scL_3^3$, $\scM_1^3$, $\scM_2^3$, $\scM_3^3$ in the variables $\bfZ=(Z_1,Z_2,Z_3)$ and $\widetilde{\bfZ}=(\widetilde{Z}_1,\widetilde{Z}_2,\widetilde{Z}_3)$, respectively, by
\begin{align*}
\scL_1^3(\bfZ) &= \xi Z_1 - Z_2, && \scM_1^3(\widetilde{\bfZ}) = \overline{\xi} \widetilde{Z}_1 - \widetilde{Z}_2, \nonumber\\
\scL_2^3(\bfZ) &= \zeta Z_2 - Z_3, && \scM_2^3(\widetilde{\bfZ}) = \overline{\zeta} \widetilde{Z}_2 - \widetilde{Z}_3, \nonumber\\
\scL_3^3(\bfZ) &= Z_2, && \scM_3^3(\widetilde{\bfZ}) = \widetilde{Z}_2.\nonumber
\end{align*}
Let $(\bfz_n)_{n\geq 1}$ be given by
\[
\forall n\in\Na \quad 
\bfz_n = (q_{s_n},q_{s_n-1},p_{s_n-1}), \quad\text{so}\quad \|\bfz_n\|_{\infty} = |q_{s_n}|.
\]
From \eqref{approx-xi} and \eqref{Ec-LemaBAD}, there is some $\kappa_6=\kappa_6(\zeta)>0$ for which
\begin{align*}
\forall n\in\clN_2 \quad \left|\scM_1^3\scM_2^3\scM_3^3(\overline{\bfz_n} )\right| &=
\left|\scL_1^3\scL_2^3\scL_3^3(\bfz_n )\right|  \nonumber\\
&= |(q_{s_n}\xi-q_{s_n-1})(\zeta q_{s_n-1} - p_{s_n-1})q_{s_n-1}| \nonumber\\
&\leq \kappa_6 \frac{|q_{s_n-1}|}{|q_{s_n-1}q_{s_n}|} 
= \frac{\kappa_6}{\| \bfz_n\|_{\infty}}. \nonumber
\end{align*}

As before, since $\|\bfz_n\|_{\infty}\to\infty$ when $n\to\infty$, Theorem \ref{ScSsThQi} assures the existence of an infinite set $\clN_3\subseteq \clN_2$ and of a vector $\bfz=(z_1,z_2,z_3)\in\Za[i]^3$, $\bfz\neq \bfzo$, satisfying
\[
\forall n\in\clN_3 \quad 
z_1q_{s_n} + z_2 q_{s_n-1} + z_3p_{s_n-1}=0.
\]
Dividing by $q_{s_n-1}$ and letting $n\to\infty$ along $\clN_3$, we get
\begin{equation}\label{EcAlfaBeta02}
\frac{1}{\xi} z_1 + z_2 + \zeta z_3 = 0 \quad\text{and}\quad z_1z_3\neq 0,
\end{equation}
because $\zeta$ and $\xi$ are irrational. Combining \eqref{EcAlfaBeta01} and \eqref{EcAlfaBeta02}, we get
\[
z_1y_2 = (z_2 + \zeta z_3)(y_1 + \zeta y_3) \quad\text{and}\quad
z_3y_3\neq 0,
\]
contradicting $[\QU(i,\zeta):\QU(i)]\geq 3$. Therefore, $\zeta$ is transcendental.

\paragraph{\S.$(w_n)_{n\geq 1}$ is strictly increasing.}

Write $t_n=w_n+u_n+v_n$ for every $n\in\Na$.
\subparagraph{\textsc{i.}} Similar to \eqref{EcCaso01Def-01}, define $\left( \zeta^{(n)} \right)_{n\geq 1}$ by
\[
\forall n\in\Na \quad \zeta^{(n)}=[0;W_nU_nV_nU_nV_nU_nV_n\ldots].
\]
Because of the hypothesis $\min_n |a_n|\geq \sqrt{8}$ and \eqref{CotaInfValidas}, the validity of the sequences defining each $\zeta^{(n)}$ is not a concern. Long but straightforward computations and \eqref{EcIdElem} tell us that each $\zeta^{(n)}$ satisfies the polynomial
\begin{align*}
P_n(X) := \begin{vmatrix} q_{w_n-1} & q_{t_n-1} \\ q_{w_n} & q_{t_n} \end{vmatrix} &X^2 - \left( \begin{vmatrix} q_{w_n-1} & p_{t_n-1} \\ q_{w_n} & p_{t_n} \end{vmatrix}  + \begin{vmatrix} p_{w_n-1} & q_{t_n-1} \\ p_{w_n} & q_{t_n} \end{vmatrix} \right)X+ \nonumber\\
&+ \begin{vmatrix} p_{w_n-1} & p_{t_n-1} \\ p_{w_n} & p_{t_n}  \end{vmatrix}.\nonumber
\end{align*}
We claim that there exists some $\kappa_1'=\kappa_1'(\zeta)>0$ such that 
\begin{equation}\label{AproxPnZeta}
\forall n\in\Na \quad 
|P_n(\zeta)| \leq \kappa_1' \frac{|q_{t_n}|}{|q_{w_n}q^2_{t_n+u_n}|}.
\end{equation}

Indeed, let $n\in\Na$. The first $t_n+u_n=w_n + 2u_n+v_n$ partial quotients of $\zeta$ and $\zeta^{(n)}$ coincide. Thus, with $\gamma$ as in Lemma \ref{PropoLemaBAD},
\[
\left| \zeta - \zeta^{(n)}\right| \leq \frac{2\gamma}{|q_{t_n+u_n}|^2},
\]
and using elementary properties of determinants we obtain
\begin{align*}
|P_n(\zeta)| &= |P_n(\zeta) - P_n(\zeta^{(n)})| \nonumber\\
&= \left| \zeta - \zeta^{(n)}\right|
\left| \begin{vmatrix} q_{w_{n}-1}\zeta-p_{w_{n}-1} & q_{t_n-1} \\ q_{w_{n}}\zeta-p_{w_{n}} & q_{t_n} \end{vmatrix}+ \begin{vmatrix} q_{w_n-1} & q_{t_n-1}\zeta^{(n)} - p_{t_n-1} \\ q_{w_n} & q_{t_n}\zeta^{(n)} - p_{t_n} \end{vmatrix} \right| \nonumber\\
&\leq 2\left| \zeta - \zeta^{(n)}\right| \left(\left| \frac{q_{t_n}}{q_{w_n}}\right| + \left|\frac{q_{w_n}}{q_{t_n}}\right|\right) \nonumber\\
&\leq \frac{2\gamma}{|q_{t_n+u_n}|^2} \left(\left| \frac{q_{t_n}}{q_{w_n}}\right| + \left|\frac{q_{w_n}}{q_{t_n}}\right|\right) \leq \kappa_1' \frac{|q_{t_n}|}{|q_{w_n}q^2_{t_n+u_n}|}.\nonumber
\end{align*}

\subparagraph{\textsc{ii.}} Let $\scL_1^1$, $\scL_2^1$, $\scL_3^1$, $\scL_4^1$, $\scM_1^1$, $\scM_2^1$, $\scM_3^1$, $\scM_4^1$ be as above and let $(\bfx_n)_{n\geq 1}$ in $\Za[i]^4$ be defined for every $n\in\Na$ by 
\[
\bfx_n:=\left(
\begin{vmatrix}
q_{w_n - 1} & q_{t_n - 1} \\
q_{w_n} & q_{t_n} \\
\end{vmatrix}, \;
\begin{vmatrix}
q_{w_n - 1} & p_{t_n - 1} \\
q_{w_n} & p_{t_n} \\
\end{vmatrix}, \;
\begin{vmatrix}
p_{w_n - 1} & q_{t_n - 1} \\
p_{w_n} & q_{t_n} \\
\end{vmatrix}, \;
\begin{vmatrix}
p_{w_n - 1} & p_{t_n - 1} \\
p_{w_n} & p_{t_n} \\
\end{vmatrix} \right),
\]
then
\begin{equation}\label{Norma_bfxn}
\forall n\in\Na \quad \|\bfx_n\|_{\infty}\leq 2 |q_{w_n}q_{t_n}|.
\end{equation}

Write $\bfx_n=(x_{n,1},x_{n,2},x_{n,3},x_{n,4})$ for $n\in\Na$. By \eqref{AproxPnZeta} and its proof, \eqref{Norma_bfxn}, Lemma \ref{LemaDaniNogueira}, and Lemma \ref{LemmaTeoGero}, there are positive constants $\kappa_2'=\kappa_2'(\zeta), \kappa_3'=\kappa_3'(\zeta)$ and an $\veps>0$ such that for all $n\in\Na$
\begin{align*}
\left| \scM^{1}_1\scM^{1}_2\scM^{1}_3\scM^{1}_4 ( \overline{\bfx_n} ) \right| &=
\left| \scL^{1}_1\scL^{1}_2\scL^{1}_3\scL^{1}_4 ( \bfx_n ) \right| \nonumber\\
&= \left| P_n(\zeta) \right| \left| \zeta x_{n,1} - x_{n,2}\right|\left| \zeta x_{n,1} - x_{n,3}\right|\left| x_{n,1}\right| 
\nonumber\\
&\leq \kappa_2' \left| \frac{q_{t_n}}{q_{w_n}q^2_{t_n+u_n}}\frac{q_{w_n}}{q_{t_n}}\frac{q_{t_n}}{q_{w_n}}q_{w_n}q_{t_n} \right| \nonumber\\
&= \kappa_2'\left| \frac{q_{t_n}^2}{q_{t_n+u_n}^2}\right| \leq \frac{\kappa_3'}{\psi^{2u_n}}\leq \frac{\kappa_3'}{|q_{w_n}q_{t_n}|^{\veps}} = \frac{\kappa_3'}{\|\bfx_n\|^{\veps}}. \nonumber
\end{align*}

Then, by Theorem \ref{ScSsThQi}, there is an $\bfx=(x_1,x_2,x_3,x_4)\in\Za[i]^4$, $\bfx\neq \bfzo$, and an infinite set $\clN_1'\subseteq \Na$ such that for all $n\in\clN_1'$
\begin{align}
0= &x_1\begin{vmatrix}
q_{w_n - 1} & q_{t_n - 1} \\
q_{w_n} & q_{t_n} \\
\end{vmatrix} + 
x_2\begin{vmatrix}
q_{w_n - 1} & p_{t_n - 1} \\
q_{w_n} & p_{t_n} \\
\end{vmatrix} + x_3\begin{vmatrix}
p_{w_n - 1} & q_{t_n - 1} \\
p_{w_n} & q_{t_n} \\
\end{vmatrix} + \nonumber\\ 
&+ x_4\begin{vmatrix}
p_{w_n - 1} & p_{t_n - 1} \\
p_{w_n} & p_{t_n} \\
\end{vmatrix}. \label{3.5}
\end{align}

Define the sequences $(Q_n)_{n\geq 1}$ and $(R_n)_{n\geq 1}$ by
\[
\forall n\in\Na \qquad
Q_n:= \frac{q_{w_n-1}q_{t_n}}{q_{w_n}q_{t_n-1}}, \quad
R_n:= \zeta - \frac{p_n}{q_n}.
\]
Note that $R_n\to 0$ as $n\to\infty$. Also, by the boundedness of $(a_n)_{n\geq 1}$ and $|a_n|\geq \sqrt{8}$ for $n\in\Na$, $(q_{w_n}/q_{w_n-1})_{n\geq 0}$ and $(q_{t_n}/q_{t_n-1})_{n\geq 0}$ are bounded and stay away from $0$. Hence, $(Q_n)_{n\geq 1}$ is bounded too. 

Divide \eqref{3.5} by $q_{w_n}q_{t_n-1}$ to obtain for every $n\in\clN_1'$
\begin{align*}
0 &= x_1(Q_n-1) + x_2\left(Q_n(\zeta - R_{t_n}) - (\zeta - R_{t_n-1})\right) + \nonumber\\
&\;+ x_3\left( Q_n(\zeta - R_{w_n-1}) - (\zeta - R_{w_n})\right)+ \nonumber\\
&\;+x_4\left( Q_n(\zeta - R_{w_n-1})(\zeta - R_{t_n}) - (\zeta - R_{w_n})(\zeta - R_{t_n-1}) \right). 
\end{align*}
Direct computations lead us to
\[
0 = (Q_n-1)(x_1+(x_2+x_3)\zeta + \zeta x_4) + \eta(n),
\]
where $\eta(n)\to 0$ when $n\to\infty$ along $\clN_1'$.

As a consequence, we obtain
\begin{equation}\label{Ec-Caso2-02}
\lim_{\substack{n\to\infty \\ n\in\clN_1'}} (Q_n-1)\left(x_1+(x_2+x_3)\zeta + \zeta^2 x_4\right) =0.
\end{equation}
Let us show that the second factor in \eqref{Ec-Caso2-02} is $0$. Take an infinite subset $\clN_1''\subseteq \clN_1'$ such that the following limits exist
\begin{equation*}
\alpha = \lim_{\substack{n\to\infty \\ n\in\clN_1''}} \frac{q_{w_n}}{q_{w_n-1}}, \quad 
\beta = \lim_{\substack{n\to\infty \\ n\in\clN_1''}} \frac{q_{t_n}}{q_{t_n-1}},
\end{equation*}
and $a=a_{w_n}\neq a_{t_n}=b$ remain constant for $n\in\clN_1''$. With $\mathfrak{K}$ as in \eqref{EcCpctDisj}, for all $n\in\clN_1''$
\[
\frac{q_{w_n}}{q_{w_n-1}} = [a;a_{w_n-1},\ldots,a_0] \in a+\mathfrak{K}, \;
\frac{q_{t_n}}{q_{t_n-1}} =[b;a_{t_n-1},\ldots,a_0]\in b+\mathfrak{K}
\]
Since $a+\mathfrak{K}$ and $b+\mathfrak{K}$ are disjoint compact sets, $\alpha\neq\beta$ and $Q_n\to \tfrac{\alpha}{\beta} \neq 1$ as $n\to\infty$ along $\clN_1''$. As a consequence, \eqref{Ec-Caso2-02} yields
\begin{equation}\label{CuadratZeta}
x_1+(x_2+x_3)\zeta + \zeta^2 x_4 =0.
\end{equation}
By $[\QU(i,\zeta):\QU(i)]\geq 3$, $x_1=x_4=x_2+x_3=0$ must hold; so \eqref{3.5} gives
\begin{equation}\label{tredive}
\forall n\in\clN_1'\quad
P_n(X) =
\begin{vmatrix} q_{w_n-1} & q_{t_n-1}\\ q_{w_n} & q_{t_n} \end{vmatrix} X^2 - 2 \begin{vmatrix} q_{w_n-1} & p_{t_n-1}\\ q_{w_n} & p_{t_n} \end{vmatrix} X + \begin{vmatrix} p_{w_n-1} & p_{t_n-1}\\ p_{w_n} & p_{t_n} \end{vmatrix}. 
\end{equation}

\subparagraph{\textsc{iii.}} Let $\scL^4_1$, $\scL^4_2$, $\scL^4_3$, $\scM^4_1$, $\scM^4_2$, $\scM^4_3$ be the linear forms in the variables $\bfY=(Y_1,Y_2,Y_3)$ and $\widetilde{\bfY}=(\widetilde{Y}_1,\widetilde{Y}_2,\widetilde{Y}_3)$ given by
\begin{align*}
\scL^4_1(\bfY) &= \zeta^2Y_1 - 2\zeta Y_2+Y_3 , \qquad
&& \scM^4_1(\widetilde{\bfY}) = \overline{\zeta}^2\widetilde{Y}_1 - 2\overline{\zeta} \widetilde{Y}_2+\widetilde{Y}_3  \nonumber\\
\scL^4_1(\bfY) &= \zeta Y_1 -Y_2 , \qquad
&& \scM^4_1(\widetilde{\bfY}) = \overline{\zeta}2\widetilde{Y}_1 - \widetilde{Y}_2  \nonumber\\
\scL^4_1(\bfY) &= Y_1, \qquad
&& \scM^4_1(\widetilde{\bfY} ) =\widetilde{Y}_1.   \nonumber
\end{align*}
Define $(\bfv_n)_{n\geq 1}$, $\bfv_n=(v_{n,1},v_{n,2},v_{n,3})\in\Za[i]^3$, by
\[
\forall n\in\Na\quad
\bfv_n:=\left( \begin{vmatrix} q_{w_n-1} & q_{t_n-1}\\ q_{w_n} & q_{t_n} \end{vmatrix}, \begin{vmatrix} q_{w_n-1} & p_{t_n-1}\\ q_{w_n} & p_{t_n} \end{vmatrix}, \begin{vmatrix} p_{w_n-1} & p_{t_n-1}\\ p_{w_n} & p_{t_n} \end{vmatrix}\right);
\]
hence,
\[
\forall n\in\Na \quad \|\bfv_n\|_{\infty}\leq 2|q_{w_n}q_{t_n}|.
\]
Using \eqref{tredive} and arguing as before, there are positive constants $\kappa_4'=\kappa_4'(\zeta)$, $\kappa_5'=\kappa_5'(\zeta), \kappa_6'=\kappa_6'(\zeta)$ such that for all $n\in\clN_1'$
\begin{align*}
\left|\scM^{4}_{1}\scM^{4}_{2}\scM^{4}_{3}(\overline{\bfv_n}) \right| &= \left|\scL^{4}_{1}\scL^{4}_{2}\scL^{4}_{3}(\bfv_n) \right| \nonumber\\
&=|P_n(\zeta) (\zeta v_{n,1}-v_{n,2}) v_{n,1}| \nonumber\\
&\leq \kappa_4'\frac{|q_{w_n}q_{t_n}|}{|q_{t_n+u_n}^2|} \leq \kappa_4'\frac{|q_{t_n}|}{|q_{t_n+u_n}|}\leq  \frac{\kappa_5'}{|q_{w_n}q_{t_n}|^{\veps}}=\frac{\kappa_6'}{\|\bfv_n\|^{\veps}}.  \nonumber
\end{align*}
Thus, Theorem \ref{ScSsThQi} ensures the existence of a vector $\bfzo\neq \bfr=(r_1, r_2, r_3)\in \Za[i]^3$ and of an infinite set $\clN_2'\subseteq \clN_1'$ such that for all $n\in\clN_2'$
\[
r_1 
\begin{vmatrix}
q_{w_n-1} & q_{t_n-1}\\
q_{w_n} & q_{t_n}
\end{vmatrix}
+ r_2
\begin{vmatrix}
q_{w_n-1} & p_{t_n-1}\\
q_{w_n} & p_{t_n}
\end{vmatrix} + r_3
\begin{vmatrix}
p_{w_n-1} & p_{t_n-1}\\
p_{w_n} & p_{t_n}
\end{vmatrix} = 0.
\]
Dividing by $q_{w_n}q_{t_n-1}$, every $n\in\clN_2'$ satisfies
\[
r_1(Q_n -1) + 
r_2 \left( Q_n\frac{p_{t_n}}{q_{t_n}} - \frac{p_{t_n-1}}{q_{t_n-1}} \right) + 
r_3 \left(Q_n \frac{p_{w_n -1}}{q_{w_n -1}} \frac{p_{t_n}}{q_{t_n}} - \frac{p_{w_n}}{q_{w_n}} \frac{p_{t_n-1}}{q_{t_n-1}} \right) = 0.
\]
As before, $(Q_n)_{n\in\clN_2'}$ has a limit point different from $1$, so 
\[
r_3\zeta^2 + r_2\zeta + r_1=0, 
\]
contradicting $[\QU(i,\zeta):\QU(i)]\geq 3$. Therefore, $\zeta$ is transcendental.
 
\section{Further results}

Other transcendence results for regular continued fractions can be translated into the HCF context with the pertinent modifications; for example, Theorem \ref{TeoBugCloverC} below can be regarded as a complex version of Theorem 1.3. in \cite{bug13autom}. 

Let $\clA$ be a finite set and $\bfa$ an infinite word on $\clA$. We say that $\bfa$ satisfies Condition $(\clubsuit)$ if it is non-periodic and there are sequences of finite words in $\clA$, $(W_n)_{n\geq 1}$, $(U_n)_{n\geq 1}$, $(V_n)_{n\geq 1}$, such that
\begin{enumerate}[i]
\item $W_nU_nV_n\widehat{U_n}$ is a prefix of $\bfa$ for every $n\in \Na$, where $\widehat{U_n}$ is the word obtained by reversing $U_n$,
\item $((|W_n|+|V_n|)/|U_n|)_{n\geq 1}$ is bounded,
\item $(|U_n|)_{n\geq 1}$ tends to infinity when $n$ does.
\end{enumerate}

\begin{teo01}\label{TeoBugCloverC}
Let $\bfa\in\Omega^{\HCF}$ satisfy Condition $(\clubsuit)$ and $|a_n|\geq \sqrt{8}$ for every $n$. Then, $\zeta=[0;a_1,a_2,\ldots]$ is transcendental.
\end{teo01}

\begin{rema1000}
\begin{enumerate}[i.]
\item As for Theorem \ref{TeoGero01} and its real counterpart, in Theorem \ref{TeoBugCloverC} we assumed that $\sean$ is bounded rather than $\sup_n |q_n|^{1/n}<+\infty$. 
\item The proof of Theorem \ref{TeoBugCloverC} requires an inequality involving the continuants. Namely, there is a constant $\kappa>0$ for which all valid prefixes $\bfa$, $\bfb$ such that $\bfa\bfb$ is a valid prefix satisfy
\[
|q_{|\bfa\bfb|-1}(\bfa\bfb)|\leq \kappa |q_{|\bfa|-1}(\bfa)q_{|\bfb|-1}(\bfb)|.
\]
The inequality follows from a classical continued fraction identity (see \cite{hensley}, Proposition 1.1). 
\item Getting rid of the condition $\min_{n\in\Na} |a_n|\geq\sqrt{8}$ in Theorem \ref{TeoBugCloverC} poses essentially the same problem as in Theorem \ref{TeoGero01}.
\end{enumerate}
\end{rema1000}
\section{Final remarks}

Our argument for showing $Q_n\not\to 1$ when $n\to\infty$ along $\clN_1'$ in the second case of Theorem \ref{TeoGero01} is simpler than the corresponding one in \cite{bug13autom}. Unfortunately, it does not shed any light upon how to omit the bounds on the partial quotients. The main complication is that, unlike regular continued fractions, for every $M>\sqrt{8}$ there are complex numbers $z=[a_0;a_1,\ldots]$ and $w=[b_0;b_1,\ldots]$ such that $|z-w|$ is arbitrarily small, $a_0\neq b_0$, $\displaystyle \liminf_{n\to\infty} |a_n|<M$, and $\displaystyle\liminf_{n\to\infty} |b_n|<M$.
 
\subsection*{Acknowledgements}
This research was supported by CONACyT (grant: 410695) and Aarhus University (St.: 2012, Pr.: 25763, Ak.: 26248). The author thanks Yann Bugeaud for the fruitful discussions, the anonymous referee whose comments lead to a significant improvement of the text, and Santiago Cabello Tueme and Juli\'an Iglesias Vargas for their feedback.


\begin{flushleft}
Gerardo Gonz\'alez Robert \par
Facultad de Ciencias, \par
Universidad Nacional Aut\'onoma de M\'exico, \par
Circuito Exterior S/N, C.U., Coyoac\'an, 04510 \par
Mexico City, Mexico \par
\texttt{gerardogonrob@ciencias.unam.mx}
\end{flushleft}


\normalsize


\begin{thebibliography}{HD82}




\baselineskip=17pt




\bibitem {adbug05} B. Adamczewski, Y. Bugeaud, \textit{On the complexity of algebraic numbers, II. Continued fractions}, Acta Math. 195 (2005), 1--20. 

\bibitem {adbulu04} B. Adamczewski, Y. Bugeaud, F. Luca, \textit{Sur la complexit\'e des nombres alg\'ebriques.} (in French) C. R. Math. Acad. Sci. Paris 339 (2004), no. 1, 11--14.

\bibitem {alshal} J.P. Allouche, J. Shallit, \textit{Automatic Sequences: Theory, Applications, Generalizations}. Cambridge. Cambridge University Press, 2003.

\bibitem {baker} A. Baker, \textit{Continued fractions of transcendental numbers.}, Mathematika 9 (1962) 1--8.

\bibitem {bosgruen} W. Bosma, D. Gruenewald, \textit{Complex numbers with bounded partial quotients}, J. Aust. Math. Soc. 93 (2012), no. 1-2, 9--20. 

\bibitem {bug13autom} Y. Bugeaud, \textit{Automatic continued fractions are transcendental or quadratic}. Ann. Sci. \'Ec. Norm. Sup\'er. (4) 46 (2013), no. 6, 1005--1022.

\bibitem {bugkim} Y. Bugeaud, D.H. Kim, \textit{A new complexity function, repetitions in Sturmian words, and irrationality exponents of Sturmian numbers.} Trans. Amer. Math. Soc. 371 (2019), no. 5, 3281--3308. 

\bibitem {dani} S. G. Dani, \textit{Continued Fractions for Complex Numbers - A General Approach}. Acta Arith. 171 (2015), no. 4, 355--369

 
\bibitem {daninog} S. G. Dani, A. Nogueira, \textit{Continued fractions for complex numbers and values of binary quadratic forms}, Trans. Amer. Math. Soc. 366 (2014), 3553--3583.


\bibitem {dodkrist} M. Dodson, S. Kristensen, \textit{Hausdorff Dimension and Diophantine Approximation}, Acta Math. 11 (2003), 187--200.

\bibitem {galois} \'E. Galois, \textit{D\'emonstration d'un th\'eor\`eme sur les fractions continues p\'eriodiques} (in French) Annales de Math\'ematiques Pures et Appliqu\'ees 19 (1829), 294--301.

\bibitem {hensley} D. Hensley, \textit{Continued Fractions}, New Jersey: World Scientific Publishing Co. Pte. Ltd., 2006.

\bibitem {hiva18} Hiary, G., Vandehey, J. \textit{Calculations of the invariant measure for Hurwitz Continued Fractions.} (2018) arXiv:1805.10151 [math.NT]


\bibitem {hines} R. Hines, \textit{Badly Approximable Numbers over Imaginary Quadratic Fields,} Acta Arith. 190(2019)101--125.


\bibitem {hur87} A. Hurwitz, \textit{\"Uber die Entwicklung complexer Gr\"ossen in Kettenbr\"uche} (in German), Acta Math. 11 (1887), 187--200.

\bibitem {juhu} J. Hurwitz, \textit{\"Uber die Reduction der Bin\"aren Quadratischen Formen mit Complexen Coefficienten und Variabeln.} (in German) 
Acta Math. 25 (1902), no. 1, 231--290. 

\bibitem {ioskraa} M. Iosifescu, C. Kraaikamp, \textit{Metrical Theory of Continued Fractions}, New Jersey: Kluwer Academic Publishers, 2002.

\bibitem {tanito} S. Ito, S. Tanaka, \textit{On a family of continued-fraction transformations and their ergodic properties.} Tokyo J. Math. 4 (1981), no. 1, 153--175.

\bibitem  {khin} A. Khinchin, A.  \textit{Continued Fractions}. New York: Dover Publications, 2006 (re-issue of the 1961 edition).

\bibitem {lakein01} R.B. Lakein, \textit{Approximation properties of some complex continued fractions}. Monatsh. Math. 77 (1973), 396--403.

\bibitem {lakein03} R.B. Lakein, \textit{Continued Fractions and Equivalent Complex Numbers}, Proc. Amer. Math. Soc. 42 (1974), 641--642.

\bibitem {leveque} W.J. LeVeque, \textit{Continued fractions and approximations in $k(i)$. I, II.} Nederl. Akad. Wetensch. Proc. Ser. A. 55 = Indagationes Math. 14, (1952). 526--535, 536--545. 

\bibitem {poitou} G. Poitou, \textit{Sur l'approximation des nombres complexes par les nombres des corps imaginaires quadratiques d\'enu\'es d'id\'eaux non principaux, particuli\`erement lorsque vaut l'algorithme d'Euclide.} (in French) Ann. Sci. Ecole Norm. Sup. (3) 70, (1953). 199--265. 

\bibitem {quef98} M. Queff\'elec, \textit{Transcendance des fractions continues de Thue-Morse.} (in French) J. Number Theory 73 (1998), no. 2, 201--211. 

\bibitem {quef00} M. Queff\'elec, \textit{Irrational numbers with automaton-generated continued fraction expansion.} Dynamical systems (Luminy-Marseille, 1998), 190--198, World Sci. Publ., River Edge, NJ, 2000.

\bibitem {aschm01} A. Schmidt, \textit{Diophantine approximation of 
complex numbers.} Acta Math. 134 (1975), 1--85

\bibitem {schm76} W.M. Schmidt, \textit{Simultaneous approximation to algebraic numbers by elements of a number field.} Monatsh. Math. 79 (1975), 55--66.

\bibitem {schm96} W.M. Schmidt, (1996), \textit{Diophantine Approximation}. Series: Lecture Notes in Mathematics 785. Berlin: Springer Verlag. 

\bibitem {tanaka} S. Tanaka, \textit{A Complex Continued Fraction Transformation and its Ergodic Properties.} Tokyo J. Math. 8 (1985), no. 1, 191--214. 
\end{thebibliography}
\end{document}